\documentclass[11pt]{article}
\usepackage{graphicx}
\usepackage{fullpage}
\usepackage{amsmath,amssymb,amsthm}
\usepackage{enumerate}
\usepackage{mathrsfs}
\usepackage{verbatim}
\usepackage{marginnote} 
\usepackage[round]{natbib}

\usepackage[T1]{fontenc}
\usepackage{microtype}

\usepackage{enumitem}
\usepackage{authblk}
\usepackage{float}

\numberwithin{equation}{section}
\newtheorem{theorem}{Theorem}
\newtheorem{lemma}{Lemma}[section]

\newtheorem{proposition}{Proposition}[section]
\newtheorem{corollary}{Corollary}[section]

\newtheorem{remark}{Remark}[section]
\newtheorem{exmp}{Example}[section]

\usepackage[usenames]{color}
\definecolor{plum}  {rgb}{.4,0,.4}
\definecolor{BrickRed} {rgb}{0.6,0,0}
\usepackage{commath}
\usepackage[normalem]{ulem}

\usepackage[plainpages=false,pdfpagelabels,colorlinks=true,linkcolor=BrickRed,citecolor=plum]{hyperref}

    \def\ddefloop#1{\ifx\ddefloop#1\else\ddef{#1}\expandafter\ddefloop\fi}

    \def\ddef#1{\expandafter\def\csname c#1\endcsname{\ensuremath{\mathcal{#1}}}}
    \ddefloop ABCDEFGHIJKLMNOPQRSTUVWXYZ\ddefloop

    \def\ddef#1{\expandafter\def\csname s#1\endcsname{\ensuremath{\mathsf{#1}}}}
    \ddefloop ABCDEFGHIJKLMNOPQRSTUVWXYZ\ddefloop

    \def\E{\mathbf{E}}

    \def\argmin{\operatornamewithlimits{arg\,min}}

    \def\bd#1{\mathbf{#1}}
    \def\bx{\bd{x}}
	\def\by{\bd{y}}
    
	\def\bB{\bd{B}}
	\def\bV{\bd{V}}

	\def\tr{\mathrm{Tr}}

\begin{document}

\title{Just Interpolate: Kernel ``Ridgeless'' Regression Can Generalize}

\author[1]{Tengyuan Liang\thanks{{\tt tengyuan.liang@chicagobooth.edu}.}}
\author[2]{Alexander Rakhlin\thanks{{\tt rakhlin@mit.edu}}}
\affil[1]{University of Chicago, Booth School of Business}
\affil[2]{Massachusetts Institute of Technology}

\date{}

\maketitle
\thispagestyle{empty}

\maketitle

\begin{abstract}
In the absence of explicit regularization, Kernel ``Ridgeless'' Regression with nonlinear kernels has the potential to fit the training data perfectly. It has been observed empirically, however, that such interpolated solutions can still generalize well on test data. We isolate a phenomenon of implicit regularization for minimum-norm interpolated solutions which is due to a combination of high dimensionality of the input data, curvature of the kernel function, and favorable geometric properties of the data such as an eigenvalue decay of the empirical covariance and kernel matrices. In addition to deriving a data-dependent upper bound on the out-of-sample error, we present experimental evidence suggesting that the phenomenon occurs in the MNIST dataset.
\end{abstract}

%
%

\section{Introduction}

According to conventional wisdom, explicit regularization should be added to the least-squares objective when the Hilbert space $\cH$ is high- or infinite-dimensional  \citep{Golub_1979,wahba1990spline,smola1998learning,shawe2004kernel,evgeniou2000regularization, de2005model, alvarez2012kernels}: 
\begin{align}
	\min_{f\in\cH} \frac{1}{n}\sum_{i=1}^n (f(x_i)-y_i)^2 + \lambda \norm{f}_{\cH}^2.
\end{align}
The regularization term is introduced to avoid ``overfitting'' since kernels provide enough flexibility to fit training data exactly (i.e. interpolate it). From the theoretical point of view, the regularization parameter $\lambda$ is a knob for balancing bias and variance, and should be chosen judiciously. Yet, as noted by a number of researchers in the last few years,\footnote{In particular, we thank M. Belkin, B. Recht, L. Rosasco, and N. Srebro for highlighting this phenomenon.} the best out-of-sample performance, empirically, is often attained by setting the regularization parameter to \emph{zero} and finding the minimum-norm solution among those that interpolate the training data. The mechanism for good out-of-sample performance of this interpolation method has been largely unclear \citep{zhang2016understanding,belkin2018understand}. 

\begin{figure}[h]
\centering
\includegraphics[width=0.6\textwidth]{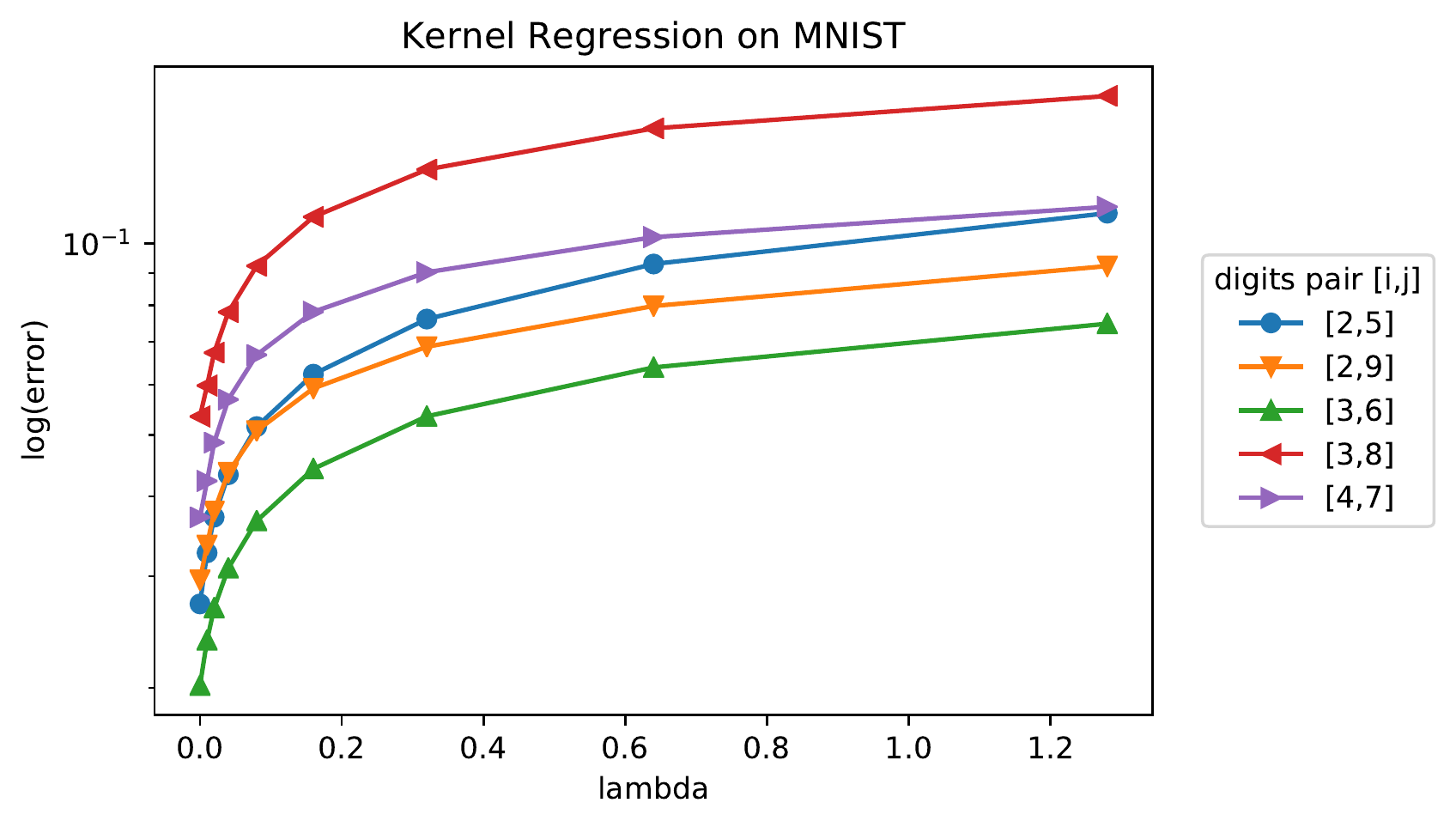}
\caption{Test performance of Kernel Ridge Regression on pairs of MNIST digits for various values of regularization parameter $\lambda$, normalized by variance of $y$ in the test set (for visualization purposes).}
\label{fig:mnist-intro}
\end{figure}

As a concrete motivating example, consider the prediction performance of Kernel Ridge Regression for various values\footnote{We take $\lambda \in \{0, 0.01, 0.02, 0.04, 0.08, 0.16, 0.32, 0.64, 1.28\}$.} of the regularization parameter $\lambda$ on subsets of the MNIST dataset. For virtually all pairs of digits, the best out-of-sample mean squared error is achieved at $\lambda=0$. Contrary to the standard bias-variance-tradeoffs picture we have in mind, the test error is monotonically decreasing as we decrease $\lambda$ (see Figure~\ref{fig:mnist-intro} and further details in Section~\ref{sec:experiments}).

We isolate what appears to be a new phenomenon of \emph{implicit regularization} for interpolated minimum-norm solutions in Kernel ``Ridgeless'' Regression. This regularization is due to the curvature of the kernel function and ``kicks in'' only for high-dimensional data and for ``favorable'' data geometry. We provide out-of-sample statistical guarantees in terms of spectral decay of the empirical kernel matrix and the empirical covariance matrix, under additional technical assumptions.  

Our analysis rests on the recent work in random matrix theory. In particular, we use a suitable adaptation of the argument of \citep{el2010spectrum} who showed that high-dimensional random kernel matrices can be approximated in spectral norm by linear kernel matrices plus a scaled identity. While the message of \citep{el2010spectrum} is often taken as ``kernels do not help in high dimensions,'' we show that such a random matrix analysis helps in explaining the good performance of interpolation in Kernel ``Ridgeless'' Regression.

\subsection{Literature Review}

Grace Wahba \citep{wahba1990spline} pioneered the study of nonparametric regression in reproducing kernel Hilbert spaces (RKHS) from the computational and statistical perspectives. One of the key aspects in that work is the role of the decay of eigenvalues of the kernel (at the population level) in rates of convergence. The analysis relies on explicit regularization (ridge parameter $\lambda$) for the bias-variance trade-off. The parameter is either chosen to reflect the knowledge of the spectral decay at the population level \citep{de2005model} (typically unknown to statistician), or by the means of cross-validation \citep{Golub_1979}. Interestingly, the explicit formula of Kernel Ridge Regression has been introduced as ``kriging'' in the literature before, and was widely used in Bayesian statistics \citep{cressie1990origins, wahba1990spline}.

In the learning theory community, Kernel Ridge Regression is known as a special case of Support Vector Regression \citep{vapnik1998statistical, shawe2004kernel, vovk2013kernel}. Notions like metric entropy \citep{cucker2002best} or ``effective dimension'' \citep{caponnetto2007optimal} were employed to analyze the guarantees on the excess loss of Kernel Ridge Regression, even when the model is misspecified. We refer the readers to \cite{gyorfi2006distribution} for more details. Again, the analysis leans crucially on the explicit regularization, as given by a careful choice of $\lambda$, for the model complexity and approximation trade-off, and mostly focusing on the fixed dimension and large sample size setting. However, to the best of our knowledge, the literature stays relatively quiet in terms of what happens to the minimum norm interpolation rules, i.e., $\lambda = 0$. As pointed out by \citep{belkin2018understand, belkin2018overfitting}, the existing bounds in nonparametric statistics and learning theory do not apply to interpolated solution either in the regression or the classification setting. In this paper, we aim to answer when and why interpolation in RKHS works, as a starting point for explaining the good empirical performance of interpolation using kernels in practice \citep{zhang2016understanding, belkin2018understand}.

\section{Preliminaries}

\subsection{Problem Formulation}

Suppose we observe $n$ i.i.d. pairs $(x_i,y_i)$, $1\leq i\leq n$, where $x_i$ are the covariates with values in a compact domain $\Omega \subset \mathbb{R}^d$ and $y_i\in \mathbb{R}$ are the responses (or, labels). Suppose the $n$ pairs are drawn from an unknown probability distribution $\mu(x,y)$. We are interested in estimating the conditional expectation function
$f_*(x) = \E(\by|\bx = x)$, which is assumed to lie in a Reproducing Kernel Hilbert Space (RKHS) $\cH$. Suppose the RKHS is endowed with the norm $\| \cdot \|_{\cH}$ and corresponding positive definite kernel $K(\cdot, \cdot): \Omega \times \Omega \rightarrow \mathbb{R}$. The interpolation estimator studied in this paper is defined as 
\begin{align}
	\label{eq:interpolation}
	\widehat{f} = \argmin_{f \in \cH} \| f \|_{\cH}, ~~ \text{s.t.}~~ f(x_i) = y_i,~\forall i \enspace.
\end{align}
Let $X \in \mathbb{R}^{n \times d}$ be the matrix with rows $x_1,\ldots,x_n$ and let $Y$ be the vector of values $y_1,\ldots,y_n$. Slightly abusing the notation, we let $K(X, X) = [K(x_i, x_j)]_{ij} \in \mathbb{R}^{n \times n}$ be the kernel matrix. Extending this definition, for $x \in \Omega$ we denote by $K(x, X) \in \mathbb{R}^{1 \times n}$ the matrix of values $[K(x,x_1),\ldots,K(x,x_n)]$. When $K(X,X)$ is invertible, solution to \eqref{eq:interpolation} can be written in the closed form:
\begin{align}
	\label{eq:interpolation_closedform}
	\widehat{f}(x) &= K(x, X) K(X, X)^{-1} Y.
\end{align}

In this paper we study the case when $K(X,X)$ is full rank, taking \eqref{eq:interpolation_closedform} as the starting point. For this interpolating estimator, we provide high-probability (with respect to a draw of $X$) upper bounds on the integrated squared risk of the form
\begin{align}
	\label{eq:target_estimation}
	\E(\widehat{f}(\bx)-f_*(\bx))^2 \leq \phi_{n,d}(X,f^*).
\end{align} 
Here the expectation is over $\bx\sim \mu$ and $Y|X$, and $\phi_{n,d}$ is a data-dependent upper bound. We remark that upper bounds of the form \eqref{eq:target_estimation} also imply prediction loss bounds for excess square loss with respect to the class $\cH$, as $\E (\widehat{f}(\bx) - f_*(\bx))^2 = \E (\widehat{f}(\bx) - \by)^2 - \E (f_*(\bx) - \by)^2$. 

\subsection{Notation and Background on RKHS}

For an operator $A$, its adjoint is denoted by $A^*$. For real matrices, the adjoint is the transpose. For any $x \in \Omega$, let $K_x: \mathbb{R} \rightarrow \cH$ be such that
\begin{align}
	\label{eq:kx}
	f(x) = \langle K_x, f \rangle_{\cH} = K_x^* f.
\end{align}
It follows that for any $x, z \in \Omega$
\begin{align}
	K(x, z) = \langle K_x, K_{z} \rangle_{\cH} = K_x^* K_{z}.
\end{align}
Let us introduce the integral operator $\mathcal{T}_\mu: L^2_\mu \rightarrow L^2_\mu$ with respect to the marginal measure $\mu(x)$:
\begin{align}
	\mathcal{T}_\mu f(z) = \int K(z, x) f(x) d\mu(x),
\end{align}
and denote the set of eigenfunctions of this integral operator by $e(x) = \{e_1(x), e_2(x), \ldots, e_p(x)\}$, where $p$ could be $\infty$. We have that
\begin{align}
	\mathcal{T}_\mu e_i = t_i e_i,~~\text{and}~~ \int e_i(x) e_j(x) d\mu(x) = \delta_{ij} \enspace.
\end{align}
Denote $T = {\rm diag}(t_1,\ldots, t_p)$ as the collection of non-negative eigenvalues. Adopting the spectral notation, 
\begin{align*}
	K(x, z) = e(x)^* T e(z). 
\end{align*}
Via this spectral characterization, the interpolation estimator \eqref{eq:interpolation} takes the following form
\begin{align}
	\widehat{f}(x) = e(x)^* T e(X) \left[ e(X)^* T e(X) \right]^{-1} Y \enspace.
\end{align}
Extending the definition of $K_x$, it is natural to define the operator $K_X: \mathbb{R}^n \rightarrow \cH$. Denote the sample version of the kernel operator to be
\begin{align}
	\widehat{\mathcal{T}} := \frac{1}{n} K_{X} K_{X}^*
\end{align}
and the associated eigenvalues to be $\lambda_j(\widehat{\mathcal{T}})$, indexed by $j$. The eigenvalues are the same as those of $\frac{1}{n} K(X, X)$. It is sometimes convenient to express $\widehat{\mathcal{T}}$ as the linear operator under the basis of eigenfunctions, in the following matrix sense
$$\widehat{T} = T^{1/2} \left( \frac{1}{n} e(X) e(X)^* \right) T^{1/2}.$$

We write $\E_{\mu}[\cdot]$ to denote the expectation with respect to the marginal $\bx\sim\mu$. Furthermore, we denote by
$$\norm{g}^2_{L^2_\mu} = \int g^2 d\mu(x) = \E_{\mu} g^2(\bx)$$
the squared $L^2$ norm with respect to the marginal distribution. The expectation $\E_{Y|X}[\cdot]$ denotes the expectation over $y_1,\ldots,y_n$ conditionally on $x_1,\ldots,x_n$.

\section{Main Result}
\label{sec:main}

We impose the following assumptions: 
\begin{enumerate}[label={\bf (A.\arabic*)}]
	\item High dimensionality: there exists universal constants $c, C \in (0, \infty)$ such that $c \leq d/n \leq C$. Denote by $\Sigma_d = \E_\mu[x_i x_i^*]$ the covariance matrix, assume that the operator norm $\| \Sigma_d \|_{\rm op} \leq 1$. 
	\item $(8+m)$-moments: $z_i := \Sigma_d^{-1/2} x_i \in \mathbb{R}^d$, $i=1,\ldots,n$, are i.i.d. random vectors. Furthermore, the entries $z_i(k), 1\leq k\leq d$ are i.i.d. from a distribution with $\E z_i(k) = 0, {\rm Var}(z_i(k)) = 1$ and $|z_i(k)| \leq C \cdot d^{\frac{2}{8+m}}$, for some $m>0$.
	\item Noise condition: there exists a $\sigma>0$ such that $\E[ (f_*(\bx)-\by)^2|\bx = x]\leq \sigma^2$ for all $x \in \Omega$.
	\item Non-linear kernel: for any $x \in \Omega$, $K(x, x) \leq M$. Furthermore, we consider the inner-product kernels of the form 
	\begin{align}
		\label{eq:inner_prod_kernel}
		K(x, x') = h\left(\frac{1}{d} \langle x, x'\rangle \right) 
	\end{align} 
	for a non-linear smooth function $h(\cdot): \mathbb{R} \rightarrow \mathbb{R}$ in a neighborhood of $0$.	
\end{enumerate}

While we state the main theorem for inner product kernels, the results follow under suitable modifications\footnote{We refer the readers to \cite{el2010spectrum} for explicit extensions to RBF kernels.} for Radial Basis Function (RBF) kernels of the form
	\begin{align}
		K(x, x')= h\left(\frac{1}{d} \| x -x' \|^2 \right). 
	\end{align}
We postpone the discussion of the assumptions until after the statement of the main theorem.

Let us first define the following quantities related to curvature of $h$:
\begin{align} 
	\label{eq:alpha-beta-gamma}
	\alpha &:= h(0) + h''(0) \frac{\tr(\Sigma_d^2)}{d^2}, \quad \beta := h'(0), \nonumber\\
	\gamma &:= h\left(\frac{\tr(\Sigma_d)}{d} \right) - h(0) - h'(0) \frac{\tr(\Sigma_d)}{d} .
\end{align}

\begin{theorem}
	\label{thm:interpolation}
	Define
	\begin{align}
		\label{eq:data-dependent-error-formula}
		\phi_{n,d}(X, f_*)  = \bV + \bB &:=  \frac{8\sigma^2 \|\Sigma_d\|_{\rm op}}{d}\sum_j \frac{\lambda_j\left(\frac{X X^* }{d} + \frac{\alpha}{\beta} 1 1^* \right)}{ \left[ \frac{\gamma}{\beta} + \lambda_j\left(\frac{X X^*}{d} +  \frac{\alpha}{\beta} 1 1^*  \right) \right]^2} \nonumber \\
		& \quad \quad \quad  +  \| f_* \|_{\cH}^2 \inf_{0\leq k \leq n}  \left\{  \frac{1}{n} \sum_{j > k} \lambda_j(K_X K_{X}^*) + 2M \sqrt{\frac{k}{n}}  \right\} \enspace.
	\end{align}
	Under the assumptions (A.1)-(A.4) and for $d$ large enough, with probability at least $1-2\delta-d^{-2}$ (with respect to a draw of design matrix $X$), the interpolation estimator \eqref{eq:interpolation_closedform}  satisfies 
	\begin{align}
		\E_{Y|X} \| \widehat{f} - f_*\|^2_{L^2_\mu} \leq \phi_{n,d}(X, f_*) + \epsilon(n, d).
	\end{align}
	Here the remainder term $\epsilon(n, d) = O(d^{-\frac{m}{8+m}} \log^{4.1} d) + O(n^{-\frac{1}{2}} \log^{0.5} (n/\delta))$. 
\end{theorem}

A few remarks are in order. First, the upper bound is \emph{data-dependent} and can serve as a certificate (assuming that an upper bound on $\sigma^2, \| f_* \|_{\cH}^2$ can be guessed) that interpolation will succeed. The bound also suggests the regimes when the interpolation method should work. The two terms in the estimate of Theorem~\ref{thm:interpolation} represent upper bounds on the variance and  bias of the interpolation estimator, respectively. Unlike the explicit regularization analysis (e.g. \citep{caponnetto2007optimal}), the two terms are not controlled by a tunable parameter $\lambda$. Rather, the choice of the non-linear kernel $K$ itself leads to an \textit{implicit control} of the two terms through curvature of the kernel function, favorable properties of the data, and high dimensionality. We remark that for the linear kernel ($h(a)=a$), we have $\gamma=0$, and the bound on the variance term can become very large in the presence of small eigenvalues. In contrast, curvature of $h$ introduces regularization through a non-zero value of $\gamma$. We also remark that the bound ``kicks in'' in the high-dimensional regime: the error term decays with both $d$ and $n$.

We refer to the favorable structure of eigenvalues of the data covariance matrix as \textit{favorable geometric properties} of the data. The first term (variance) is small when the data matrix enjoys certain decay of the eigenvalues, thanks to the implicit regularization $\gamma$. The second term (bias) is small when the eigenvalues of the kernel matrix decay fast or the kernel matrix is effectively low rank. Note that the quantities $\alpha, \beta$ are constants, and $\gamma$ scales with $(\tr(\Sigma_d)/d)^2$. We will provide a detailed discussion on the trade-off between the bias and variance terms for concrete examples in Section~\ref{sec:data-dependent-bound}.

We left the upper bound of Theorem~\ref{thm:interpolation} in a data-dependent form for two reasons. First, an explicit dependence on the data tells us whether interpolation can be statistically sound on the given dataset. Second, for general spectral decay, current random matrix theory falls short of characterizing the spectral density non-asymptotically except for special cases \citep{bose2003limiting, el2010spectrum}.

\paragraph{Discussion of the assumptions} 
\begin{itemize}
	\item The assumption in (A.1) that $c\leq d/n \leq C$ emphasizes that we work in a high-dimensional regime where $d$ scales on the order of $n$. This assumption is used in the proof of \citep{el2010spectrum}, and the particular dependence on $c,C$ can be traced in that work if desired. Rather than doing so, we ``folded'' these constants into mild additional power of $\log d$. The same goes for the assumption on the scaling of the trace of the population covariance matrix.
	\item The assumption in (A.2) that $Z_i(k)$ are i.i.d. across $k=1,\ldots,d$ is a strong assumption that is required to ensure the favorable high-dimensional effect. Relaxing this assumption is left for future work. 
	\item The existence of $(8+m)$-moments for $|z_i(k)|$ is enough to ensure $|z_i(k)| \leq C \cdot d^{\frac{2}{8+m}}$ for $1\leq i\leq n, 1\leq k\leq d$ almost surely (see, Lemma 2.2 in \cite{yin1988limit}). Remark that the assumption of existence of $(8+m)$-moments in (A.2) is relatively weak. In particular, for bounded or subgaussian variables, $m=\infty$ and the error term $\epsilon(n,d)$ scales as $d^{-1} + n^{-1/2}$, up to log factors. See Lemma~\ref{lem:1byn-gaussian} for an explicit calculation in the Gaussian case. 
	\item Finally, as already mentioned, the main result is stated for the inner product kernel, but can be extended to the RBF kernel using an adaptation of the analysis in \citep{el2010spectrum}.
\end{itemize}

\section{Behavior of the Data-dependent Bound}
\label{sec:data-dependent-bound}

In this section, we estimate, both numerically and theoretically, the non-asymptotic data-dependent upper bound in Theorem~\ref{thm:interpolation} in several regimes. To illustrate the various trade-offs, we divide the discussion into two main regimes: $n>d$ and $n<d$. Without loss of generality, we take as an illustration the non-linearity $h(t) = \exp(2t)$ and $K(x, x') = \exp(2\langle x, x'\rangle/d)$, with the implicit regularization ${\bf r}:= \gamma/\beta \asymp \left(\tr(\Sigma_d)/d \right)^2$. In our discussion, we take both $n$ and $d$ large enough so that the residuals in Theorem~\ref{thm:interpolation} are negligible. The main theoretical results in this section, Corollaries~\ref{coro:general-n-d} and \ref{coro:general-d-n}, are direct consequences of the data-dependent bound in Theorem~\ref{thm:interpolation}.

\paragraph{Case $n>d$}
We can further bound the variance and the bias, with the choice $k=0$, as
\begin{align}
	\label{eq:V_n_d}
	{\bf V} &\precsim \frac{1}{d} \sum_j \frac{\lambda_j\left(\frac{X X^* }{d}  \right)}{ \left[ {\bf r}  + \lambda_j\left(\frac{X X^*}{d}  \right) \right]^2} = \frac{1}{n} \sum_{j=1}^d \frac{\lambda_j\left(\frac{X X^* }{n}  \right)}{ \left[ \frac{d}{n} {\bf r}   + \lambda_j\left(\frac{X X^*}{n}  \right) \right]^2}, \\
	\label{eq:B_n_d}
	{\bf B} &\precsim \frac{1}{n} \sum_{j=1}^n \lambda_j(K_X K_{X}^*) \asymp {\bf r} + \frac{1}{d} \sum_{j=1}^d \lambda_j\left(\frac{X X^* }{n}\right).
\end{align}

We first illustrate numerically the bias-variance trade-off by varying the geometric properties of the data in terms of the population spectral decay of $\bx$. We shall parametrize the eigenvalues of the covariance, for $0<\kappa<\infty$, as
\begin{align*}
	\lambda_j(\Sigma_d) = \left(1 - ((j-1)/d)^{\kappa}\right)^{1/\kappa},  1\leq j \leq d.
\end{align*}
The parameter $\kappa$ controls approximate ``low-rankness'' of the data: the closer $\kappa$ is to $0$, the faster does the spectrum of the data decay. This is illustrated in the top row of Figure~\ref{fig:n_m_d} on page \pageref{fig:n_m_d}. By letting $\kappa \rightarrow 0$, ${\bf r}$ can be arbitrary small, as
\begin{align*}
	\frac{\tr(\Sigma_d)}{d} \asymp \int_0^1 (1 - t^{\kappa})^{1/\kappa} dt  = \frac{\Gamma(1+1/\kappa)^2}{\Gamma(1+2/\kappa)} \in [0,1] .
\end{align*}
We will focus on three cases, $\kappa \in \{e^{-1}, e^{0}, e^{1}\}$, for the decay parameter, and values $d = 100$, $n \in \{500, 2000\}$. The data-dependent upper bounds on $\bV$ and $\bB$ are summarized in Table~\ref{tab:simuls_n_d}. More detailed plots are postponed to Figure~\ref{fig:n_m_d} (in this figure, we plot the ordered eigenvalues and the spectral density for both the population and empirical covariances). Table~\ref{tab:simuls_n_d} shows that as $\kappa$ increases (a slower spectral decay), the implicit regularization parameter becomes larger, resulting in a decreasing variance and an increasing bias.

\begin{figure}[ht]
\centering

\includegraphics[width=0.5\textwidth]{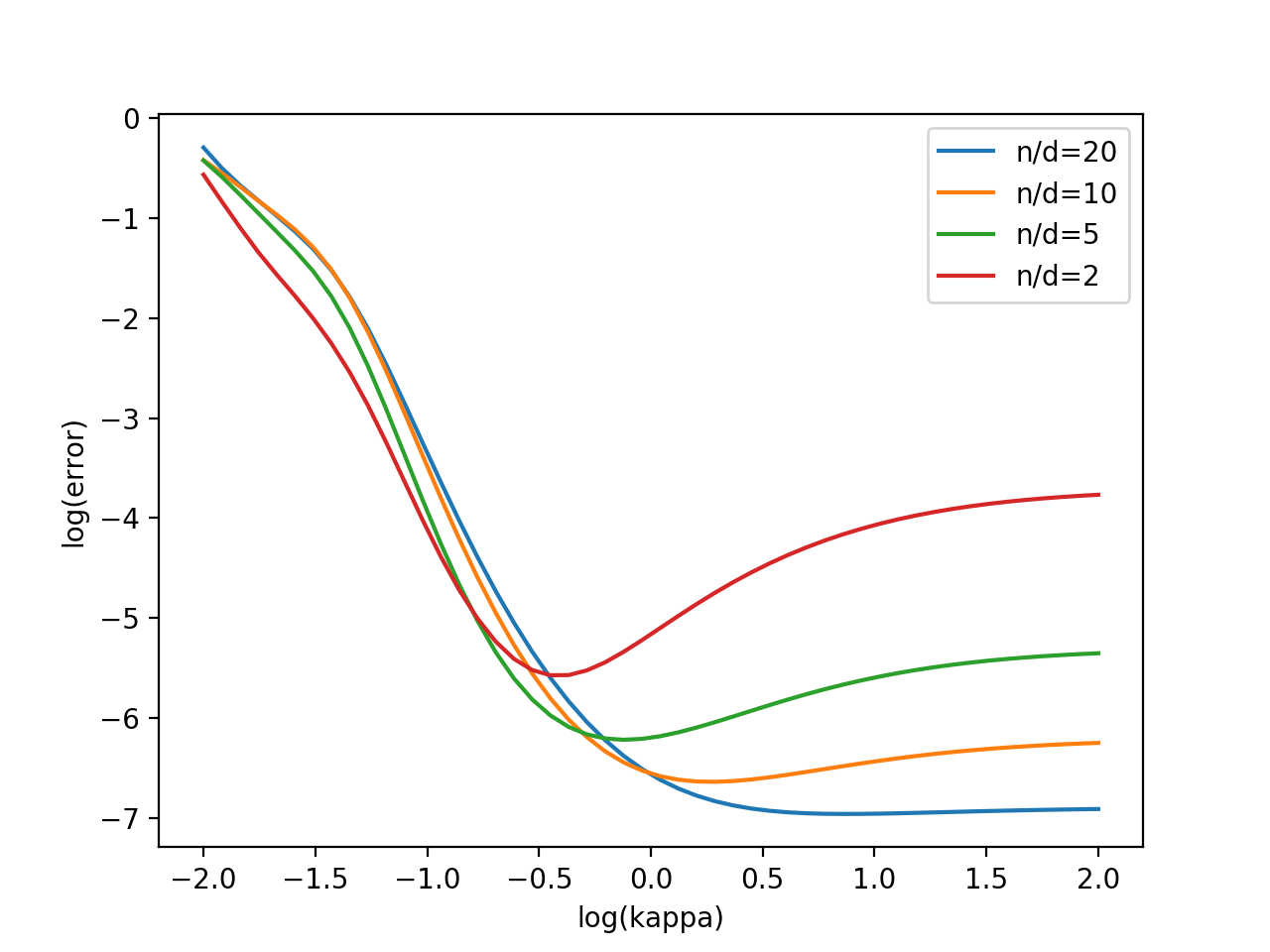}

\caption{Generalization error as a function of varying spectral decay. Here $d = 200$, $n = 400, 1000, 2000, 4000$. }
\label{fig:gen_error_n_d}
\end{figure}

\begin{figure}[ht]
\centering

\includegraphics[width=0.5\textwidth]{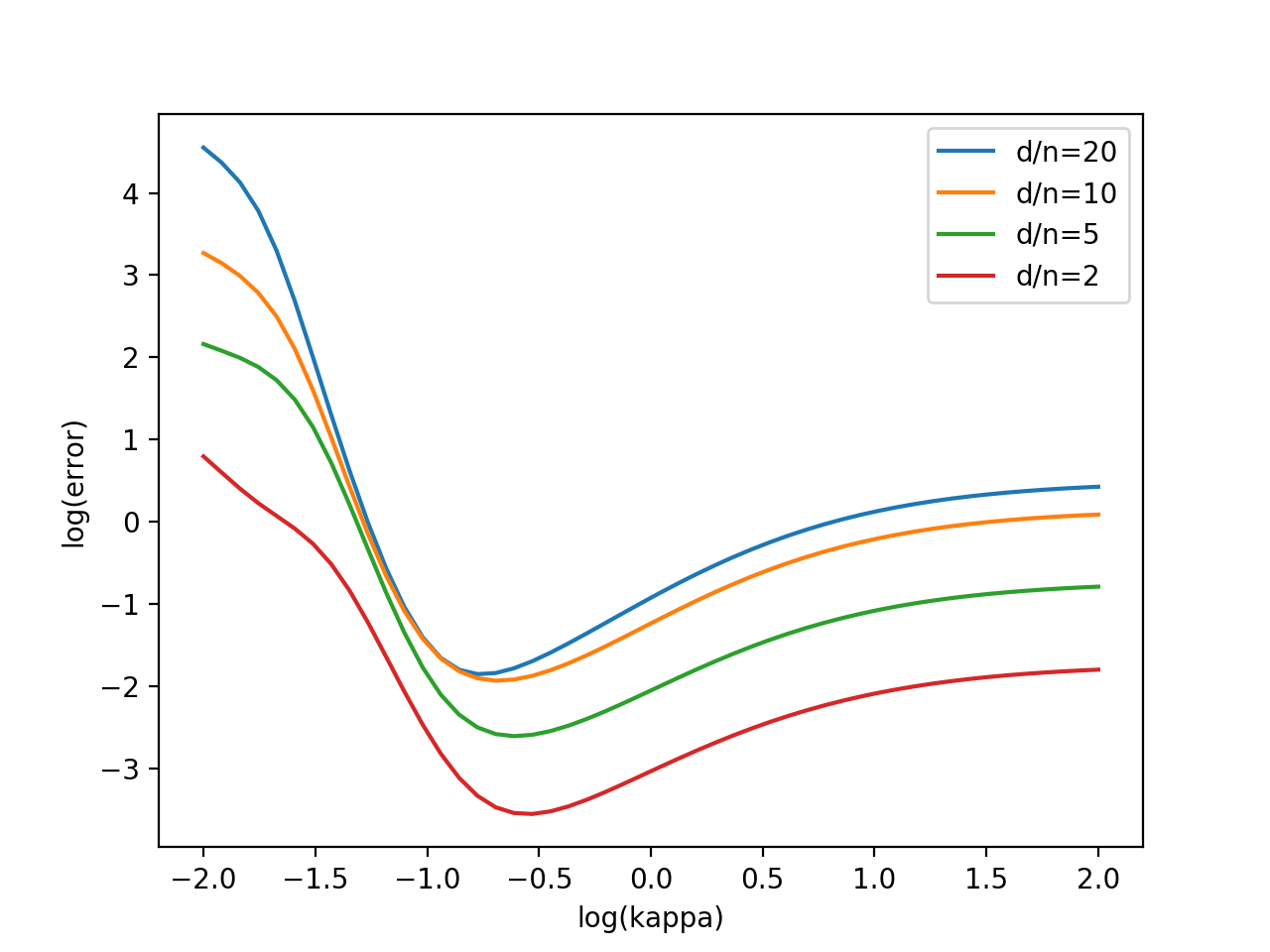}

\caption{Generalization error as a function of varying spectral decay. Here $n=200$, $d = 400, 1000, 2000, 4000$.}
\label{fig:gen_error_d_n}
\end{figure}

We also perform simulations to demonstrate the trade-off between bias and variance in the generalization error. The result is shown in Figure~\ref{fig:gen_error_n_d}.
For each choice of $(n,d)$ pair, we vary the spectral decay of the kernel by changing gradually $\kappa \in [e^{-2}, e^{2}]$, and plot the generalization error on the log scale.
We postpone the experiment details to Section~\ref{sec:synthetic}, but point out important phenomenona observed in Figures~\ref{fig:gen_error_n_d}-\ref{fig:gen_error_d_n}: (1) an extremely fast spectral decay (small $\kappa$) will generate insufficient implicit regularization that would hurt the generalization performance due to a large variance term; (2) a very slow spectral decay (large $\kappa$) will result in a large bias, which can also hurt the generalization performance; (3) certain favorable spectral decay achieves the best trade-off, resulting in the best generalization error.

\begin{table}[ht]
	\centering
	\caption{}\label{tab:simuls_n_d}
	\small
	\renewcommand{\arraystretch}{1.2}
	\begin{tabular}{lcrrrrr}
		\multicolumn{7}{c}{\normalsize Case $n > d$: variance bound $\bf V$ \eqref{eq:V_n_d}, bias bound $\bf B$ \eqref{eq:B_n_d}} \\
		\hline
		\hline
		& & \multicolumn{2}{c}{$n/d=5$} &	&			\multicolumn{2}{c}{$n/d=20$}       \\		
		\cline{3-4} \cline{6-7}
		Spectral Decay &  Implicit Reg &	$\bf V$     &     	$\bf B$        &&    		$\bf V$      &     	$\bf B$        \\     
		\hline															
		$\kappa=e^{-1}$ &  0.005418    &     14.2864    &     	0.07898         &&     		 9.4980    &     	0.07891         \\     
		$\kappa=e^{0}$  &  0.2525      &     0.4496     &     	0.7535          &&     	    0.1748	   &     	0.7538         \\ 
		$\kappa=e^{1}$  &  0.7501      &     0.1868     &     	1.6167          &&     		0.05835    &     	1.6165         \\     
		\hline
		\hline				
	\end{tabular}
\end{table}

We now theoretically demonstrate scalings within the $n > d$ regime when both $\bf V$ and $\bf B$ vanish. For simplicity, we consider Gaussian $\bx$.

\begin{corollary}[General spectral decay: $n>d$]
	\label{coro:general-n-d}
	Consider general eigenvalue decay with $\| \Sigma_d \|_{\rm op}\leq 1$. Then with high probability,
	\begin{align*}
		{\bf V} \precsim \frac{\tr(\Sigma_d^{-1})}{n}, \quad {\bf B} \precsim {\bf r} + \frac{\tr(\Sigma_d)}{d}.
	\end{align*}
\end{corollary}

To illustrate the behavior of the estimates in Corollary~\ref{coro:general-n-d}, consider the following assumptions on the population covariance matrix:

\begin{exmp}[Low rank]
	\label{exmp:low-rank}
	Let $\Sigma_d = {\rm diag}(1, \ldots, 1, 0, \ldots, 0)$ with $\epsilon d$ ones, $\epsilon\in(0,1)$. In this case
	${\bf r} = \epsilon^2$, and $\lambda_j(XX^*/n) \geq (1 - \sqrt{\epsilon d/n})^2$ with high probability by standard results in random matrix theory. Then
	\begin{align*}
		{\bf V} \precsim \frac{\epsilon d}{n} \frac{ (1 - \sqrt{\epsilon d/n})^2}{ \left( \epsilon^2 d/n +  (1 - \sqrt{\epsilon d/n})^2 \right)^2} \asymp \frac{d}{n}\epsilon , \quad {\bf B} \precsim \epsilon^2 + \epsilon.
	\end{align*}
	Therefore, as $\epsilon \rightarrow 0$, both terms vanish for $n > d$.  
\end{exmp}

\begin{exmp}[Approx. low rank]
	\label{exmp:approx-low-rank}
	Let $\Sigma_d = {\rm diag}(1, \epsilon,  \ldots, \epsilon)$ for small $\epsilon>0$. In this case,
	${\bf r} = \epsilon^2$ and $\lambda_j(XX^*/n) \geq \epsilon (1 - \sqrt{d/n})^2$ with high probability. Then
	\begin{align*}
		{\bf V} \precsim \frac{d}{n} \frac{\epsilon (1 - \sqrt{d/n})^2}{ \left( \epsilon^2 d/n + \epsilon (1 - \sqrt{d/n})^2 \right)^2} \asymp \frac{d}{n} \frac{1}{\epsilon}, \quad {\bf B} \precsim \epsilon^2 + \epsilon.
	\end{align*}
	For instance, for $\epsilon \asymp (d/n)^{1/2}$, both terms vanish for $n \gg d$.  
\end{exmp}

\begin{exmp}[Nonparametric slow decay]
	\label{exmp:nonparam}
	Consider $\lambda_j(\Sigma_d) = j^{-\alpha}$ for $0<\alpha<1$. Then ${\bf r} \asymp d^{-2\alpha}$. One can bound w.h.p. (see \eqref{pf:example-nonparam})
	\begin{align*}
		{\bf V} 
		&\asymp \frac{1}{n} \int_0^d t^{\alpha} dt \asymp \frac{d^{\alpha+1}}{n}, \quad {\bf B} \precsim d^{-2\alpha} + d^{-\alpha}.
	\end{align*}
	Balancing the two terms, one obtains a nonparametric upper bound
	$
		n^{-\frac{\alpha}{2\alpha+1}}.
	$
	A similar analysis can be carried out for $\alpha\geq 1$.
\end{exmp}

\paragraph{Case $d>n$}

In this case, we can further bound the variance and the bias, with the choice $k=0$, as
\begin{align}
	\label{eq:V_d_n}
	{\bf V} &\precsim \frac{1}{d} \sum_{j=1}^n \frac{\lambda_j\left(\frac{X X^* }{d}  \right)}{ \left[ {\bf r}  + \lambda_j\left(\frac{X X^*}{d}  \right) \right]^2}, \\
	\label{eq:B_d_n}
	{\bf B} &\precsim \frac{1}{n} \sum_{j=1}^n \lambda_j(K_X K_{X}^*) \asymp {\bf r} + \frac{1}{n} \sum_{j=1}^n \lambda_j\left(\frac{X X^* }{d}\right).
\end{align}

We first numerically illustrate the trade-off between the variance and the bias upper bounds. We consider three cases $\kappa \in \{ e^{-1}, e^{0}, e^{1} \}$, and $d = 2000$, $n \in \{400, 100\}$. As before, we find a trade-off between $\bV$ and $\bB$ with varying $\kappa$; the results are summarized in Table~\ref{tab:simuls_d_n}. Additionally, Figure~\ref{fig:d_m_n} provides a plot of the ordered eigenvalues, as well as spectral density for both the population and empirical covariances. As one can see, for a general eigenvalue decay, the spectral density of the population and the empirical covariance can be quite distinct. We again plot the generalization error in Figure~\ref{fig:gen_error_d_n} as a function of $\kappa$.

\begin{table}[ht]
	\centering
	\caption{}\label{tab:simuls_d_n}
	\small
	\renewcommand{\arraystretch}{1.2}
	\begin{tabular}{lcrrrrr}
		\multicolumn{7}{c}{\normalsize Case $d > n$: variance bound $\bf V$ \eqref{eq:V_d_n}, bias bound $\bf B$ \eqref{eq:B_d_n}} \\
		\hline
		\hline
		& & \multicolumn{2}{c}{$d/n=5$} &	&			\multicolumn{2}{c}{$d/n=20$}       \\		
		\cline{3-4} \cline{6-7}
		Spectral Decay &  Implicit Reg &	$\bf V$     &     	$\bf B$        &&    		$\bf V$      &     	$\bf B$        \\     
		\hline															
		$\kappa=e^{-1}$ &  0.005028    &     3.9801    &     	0.07603        &&       	0.7073	    &     	0.07591         \\     
		$\kappa=e^{0}$  &  0.2503      &     0.1746    &     	0.7513         &&           0.04438		&     	0.7502         \\ 
		$\kappa=e^{1}$  &  0.7466      &     0.06329   &     	1.6106         &&     	    0.01646	    &     	1.6102         \\     
		\hline
		\hline
	\end{tabular}
\end{table}

We now theoretically showcase an example in the $d \gg n$ regime where both $\bf V$ and $\bf B$ vanish. Again consider $\bx$ being Gaussian for simplicity.
\begin{corollary}[General spectral decay: $d > n$]
	\label{coro:general-d-n}
	With high probability, it holds that
	\begin{align*}
		{\bf V} \precsim \frac{n}{d} \frac{1}{4 {\bf r}}, \quad  
		{\bf B} \precsim  {\bf r} +  \frac{\tr(\Sigma_d)}{d}.
	\end{align*}
	The variance bound follows from the fact that $\frac{t}{({\bf r} + t)^2} \leq \frac{1}{4{\bf r}}$ for all $t$.
\end{corollary}

\begin{exmp}[Favorable spectral decay for $d\gg n$]
	\label{exmp:favor-decay}
	Recall $\tr(\Sigma_d)/d = {\bf r}^{1/2}$. With the choice ${\bf r} = (n/d)^{2/3}$, both terms vanish for $d \gg n$ as 
	\begin{align*}
		{\bf V} &\precsim \frac{n}{d} \frac{1}{4 {\bf r}}, \quad {\bf B} \precsim {\bf r}^{1/2}.
	\end{align*}
	In this case, the spectrum satisfies $\tr(\Sigma_d)/d = O((n/d)^{1/3})$.
\end{exmp}

\section{Proofs}

To prove Theorem~\ref{thm:interpolation}, we decompose the mean square error into the bias and variance terms (Lemma~\ref{lem:decomposition}), and provide data-dependent bound for each (Sections \ref{sec:variance} and \ref{sec:bias}). 

\subsection{Bias-Variance Decomposition}

The following is a standard bias-variance decomposition for an estimator. We remark that it is an equality, and both terms have to be small to ensure the desired convergence.

\begin{lemma}
	\label{lem:decomposition}
	The following decomposition for the interpolation estimator \eqref{eq:interpolation_closedform} holds
	\begin{align}
		\E_{Y|X} \| \widehat{f} - f_*\|^2_{L^2_\mu} = \bV + \bB , 
	\end{align}
	where
	\begin{align}
	\bV&:= \int \E_{Y|X} \left| K_x^* K_{X} (K_{X}^* K_{X})^{-1} (Y - \E[Y|X]) \right|^2 d\mu(x), \label{eq:variance}\\
	\bB&:= \int \left| K_x^* \left[ K_{X} (K_{X}^* K_{X})^{-1} K_{X}^* - I \right] f_* \right|^2 d\mu(x). 
	\label{eq:bias}
	\end{align}
\end{lemma}
 
\begin{proof}[Proof of Lemma~\ref{lem:decomposition}]
	Recall the closed form solution of the interpolation estimator: 
	\begin{align*}
		\widehat{f}(x) &= K_x^* K_{X} (K_{X}^* K_{X})^{-1} Y = K(x, X) K(X, X)^{-1} Y.
	\end{align*}
	Define $E = Y - \E[Y|X] = Y - f_*(X)$. Since $\E_{Y|X} E = 0$, we have
	\begin{align*}
		\widehat{f}(x) - f_*(x) &=  K_x^* K_{X} (K_{X}^* K_{X})^{-1} E + K_x^* \left[ K_{X} (K_{X}^* K_{X})^{-1} K_{X}^* - I \right] f_* \\
		\E_{Y|X} (\widehat{f}(x) - f_*(x))^2 &= \E_{Y|X} \left(K_x^* K_{X} (K_{X}^* K_{X})^{-1} E \right)^2 + \left| K_x^* \left[ K_{X} (K_{X}^* K_{X})^{-1} K_{X}^* - I \right] f_* \right|^2.
	\end{align*}
	Using Fubini's Theorem, 
	\begin{align*}
	 	&\E_{Y|X} \| \widehat{f} - f_*\|^2_{L^2_\mu}  =  \int \E_{Y|X} (\widehat{f}(x) - f_*(x))^2 d\mu(x) \\
		& = \int \E_{Y|X} \left| K_x^* K_{X} (K_{X}^* K_{X})^{-1} E \right|^2 d\mu(x) + \int \left| K_x^* \left[ K_{X} (K_{X}^* K_{X})^{-1} K_{X}^* - I \right] f_* \right|^2 d\mu(x). 
	\end{align*}
\end{proof}

\subsection{Variance}
\label{sec:variance}

In this section, we provide upper estimates on the variance part $\bV$ in \eqref{eq:variance}. 
\begin{theorem}[Variance]
	\label{thm:variance}
	Let $\delta\in (0,1)$. Under the assumptions (A.1)-(A.4), with probability at least $1 - \delta - d^{-2}$ with respect to a draw of $X$, 
	\begin{align}
		\label{eq:var-formula}
		\bV \leq \frac{8\sigma^2 \| \Sigma_d \|}{d} \sum_j \frac{\lambda_j\left(\frac{X X^*}{d} + \frac{\alpha}{\beta} 11^* \right)}{ \left[ \frac{\gamma}{\beta} + \lambda_j\left(\frac{X X^*}{d} + \frac{\alpha}{\beta} 11^*\right) \right]^2} + \frac{8\sigma^2}{\gamma^2} d^{-(4\theta - 1)} \log^{4.1} d,
	\end{align}
	for $\theta = \frac{1}{2} - \frac{2}{8+m}$ and for $d$ large enough.
\end{theorem}

\begin{remark}
	Let us discuss the first term in Eq.~\eqref{eq:var-formula} and its role in  implicit regularization induced by the curvature of the kernel, eigenvalue decay, and high dimensionality. In practice, the data matrix $X$ is typically centered, so $1^*X = 0$. Therefore the first term is effectively
	\begin{align*}
		\sum_{j} f_r\left( \lambda_j\left(\frac{X X^*}{d} \right) \right), ~\text{where}~ f_r(t) := \frac{t}{(r+t)^2} \leq \frac{1}{4r}.
	\end{align*}
	This formula explains the effect of implicit regularization, and captures the ``effective rank'' of the training data $X$. We would like to emphasize that this measure of complexity is distinct from the classical notion of effective rank for regularized kernel regression \citep{caponnetto2007optimal}, where the ``effective rank'' takes the form $\sum_j g_r(t_j)$ with $g_r(t) = t/(r+t)$, with $t_j$ is the eigenvalue of the population integral operator $\mathcal{T}$.
\end{remark}

\begin{proof}[Proof of Theorem~\ref{thm:variance}]
	From the definition of $\bV$ and $E[Y|X] = f_*(X)$,
	\begin{align*}
		\bV &=  \int \E_{Y|X} \tr\left(  K_x^* K_{X} (K_{X}^* K_{X})^{-1} (Y - f_*(X)) (Y -  f_*(X))^*  (K_{X}^* K_{X})^{-1} K_{X}^* K_x \right) d\mu(x) \\
		& \leq \int \| (K_{X}^* K_{X})^{-1} K_{X}^* K_x \|^2 \| \E_{Y|X} \left[(Y - f_*(X)) (Y -  f_*(X))^* \right]\| d\mu(x).
	\end{align*}
	Due to the fact that $\E_{Y|X}\left[ (Y_i - f_*(X_i))(Y_j -f_*(X_j)) \right] = 0$ for $i \neq j$, and $\E_{Y|X}\left[ (Y_i - f_*(X_i))^2 \right] \leq \sigma^2$, we have that $\| \E_{Y|X} \left[(Y - f_*(X)) (Y -  f_*(X))^* \right]\| \leq \sigma^2$ and thus
	\begin{align*}
		\bV &\leq \sigma^2 \int \| (K_{X}^* K_{X})^{-1} K_{X}^* K_x \|^2  d\mu(x) = \sigma^2 \E_\mu \| K(X, X)^{-1} K(X, \bx) \|^2.
	\end{align*}
	
	Let us introduce two quantities for the ease of derivation. For $\alpha, \beta, \gamma$ defined in \eqref{eq:alpha-beta-gamma}, let
	\begin{align}
		K^{\rm lin}(X, X) &:=  \gamma I + \alpha 11^T + \beta \frac{X X^*}{d} \in \mathbb{R}^{n \times n}, \\
		K^{\rm lin}(X, x) &:=  \beta \frac{X x^*}{d} \in \mathbb{R}^{n \times 1}, 
	\end{align}
	and $K^{\rm lin}(x, X)$ being the transpose of $K^{\rm lin}(X, x)$.
	By Proposition \ref{prop:Karoui-result}, with probability at least $1-\delta-d^{-2}$, for $\theta = \frac{1}{2} - \frac{2}{8+m}$ the following holds
	\begin{align*}
		\left\| K(X, X) - K^{\rm lin}(X, X) \right\| \leq d^{-\theta} (\delta^{-1/2} + \log^{0.51} d).
	\end{align*}
	As a direct consequence, one can see that
	\begin{align}
		\left\| K(X, X)^{-1} \right\| &\leq  \frac{1}{\gamma - d^{-\theta} (\delta^{-1/2} + \log^{0.51} d)} \leq \frac{2}{\gamma}, \label{eq:nbyn1}\\
		\left\| K(X, X)^{-1} K^{\rm lin}(X, X) \right\| &\leq 1 + \| K(X, X)^{-1} \| \cdot \| K(X, X) -  K^{\rm lin}(X, X) \| \nonumber \\ 
		& \leq \frac{\gamma}{\gamma - d^{-\theta} (\delta^{-1/2} + \log^{0.51} d)} \leq 2, \label{eq:nbyn2}
	\end{align}
	provided $d$ is large enough, in the sense that
	$$
	d^{-\theta} (\delta^{-1/2} + \log^{0.51} d) \leq \gamma/2.
	$$
	By Lemma~\ref{lemma:1byn} (for Gaussian case, Lemma~\ref{lem:1byn-gaussian}), 
	\begin{align}
		\label{eq:1byn}
		\E_\mu \left\| K(\bx, X) -  K^{\rm lin}(\bx, X) \right\|^2 \leq d^{-(4 \theta - 1)} \log^{4.1} d.
	\end{align}

	Let us proceed with the bound
	\begin{align*}
		\bV & \leq  \sigma^2 \E_\mu \| K(X, X)^{-1} K(X, \bx) \|^2  \\
		& \leq 2\sigma^2 \E_\mu \| K(X, X)^{-1} K^{\rm lin}(X, \bx) \|^2 + 2 \sigma^2 \left\| K(X, X)^{-1} \right\|^2 \cdot \E_\mu \|K(X, \bx) -  K^{\rm lin}(X, \bx)\|^2 \\
		& \leq 2\sigma^2 \left\| K(X, X)^{-1} K^{\rm lin}(X, X) \right\|^2 \E_\mu \| K^{\rm lin}(X, X)^{-1} K^{\rm lin}(X, \bx) \|^2 + \frac{8\sigma^2}{\gamma^2} d^{-(4\theta - 1)} \log^{4.1} d \\
		& \leq 8\sigma^2  \E_\mu \| K^{\rm lin}(X, X)^{-1} K^{\rm lin}(X, \bx) \|^2 +  \frac{8\sigma^2}{\gamma^2} d^{-(4\theta - 1)} \log^{4.1} d
	\end{align*}
	where the the third inequality relies on \eqref{eq:1byn} and \eqref{eq:nbyn1}, and the fourth inequality follows from \eqref{eq:nbyn2}.
	
	One can further show that
	\begin{align*}
	 &\E_\mu \| K^{\rm lin}(X, X)^{-1} K^{\rm lin}(X, \bx) \|^2 \\
	 &= \E_{\mu}  \tr\left( \left[\gamma I + \alpha 1 1^* + \beta \frac{XX^*}{d} \right]^{-1} \beta \frac{X \bx}{d}  \beta \frac{\bx^{*} X^{*}}{d} \left[\gamma I + \alpha 1 1^* + \beta \frac{XX^*}{d} \right]^{-1} \right) \\
	 	&= \tr\left( \left[\gamma I + \alpha 1 1^* + \beta \frac{XX^*}{d} \right]^{-1} \beta^2 \frac{X \Sigma_d X^{*}}{d^2}  \left[\gamma I + \alpha 1 1^* + \beta \frac{XX^*}{d} \right]^{-1} \right)\\
		& \leq \frac{1}{d} \| \Sigma_d \| \tr\left( \left[\gamma I + \alpha 1 1^* + \beta \frac{X^* X}{d} \right]^{-1} \beta^2 \frac{X^*X}{d} \left[\gamma I + \alpha 1 1^* + \beta \frac{X^* X}{d} \right]^{-1}  \right) \\
		& \leq \frac{1}{d} \| \Sigma_d \| \tr\left( \left[\gamma I + \alpha 1 1^* + \beta \frac{X^* X}{d} \right]^{-1} \left[\beta^2 \frac{X^*X}{d} +  \alpha \beta 1 1^*  \right] \left[\gamma I + \alpha 1 1^* + \beta \frac{X^* X}{d} \right]^{-1}  \right) \\
		& = \frac{1}{d} \| \Sigma_d \|  \sum_j \frac{\lambda_j\left(\frac{X X^*}{d} + \frac{\alpha}{\beta} 11^* \right)}{ \left[ \frac{\gamma}{\beta} + \lambda_j\left(\frac{X X^*}{d} + \frac{\alpha}{\beta} 11^*\right) \right]^2}.
	\end{align*}
	
	We conclude that with probability at least $1-\delta-d^{-2}$,
	\begin{align}
		\bV &\leq  8\sigma^2  \E_\mu \| K^{\rm lin}(X, X)^{-1} K^{\rm lin}(X, \bx) \|^2 +  \frac{8\sigma^2}{\gamma^2} d^{-(4\theta - 1)} \log^{4.1} d \\
		& \leq \frac{8\sigma^2 \| \Sigma_d \|}{d} \sum_j \frac{\lambda_j\left(\frac{X X^*}{d} + \frac{\alpha}{\beta} 11^* \right)}{ \left[ \frac{\gamma}{\beta} + \lambda_j\left(\frac{X X^*}{d} + \frac{\alpha}{\beta} 11^*\right) \right]^2} + \frac{8\sigma^2}{\gamma^2} d^{-(4\theta - 1)}\log^{4.1} d
	\end{align}
	for $d$ large enough. 
\end{proof}

\subsection{Bias}
\label{sec:bias}


\begin{theorem}[Bias]
	\label{thm:bias}
	Let $\delta\in(0,1)$. The bias, under the only assumptions that $K(x, x)\leq M$ for $x \in \Omega$, and $X_i$'s are i.i.d. random vectors, is upper bounded as
	\begin{align}
			\bB \leq \| f_* \|_{\cH}^2 \cdot \inf_{0\leq k \leq n}  \left\{  \frac{1}{n} \sum_{j > k} \lambda_j(K(X, X)) + 2 \sqrt{\frac{k}{n}} \sqrt{ \frac{\sum_{i=1}^n K(x_i, x_i)^2}{n}}  \right\} + 3M \sqrt{\frac{\log 2n/\delta}{2n}}, 
	\end{align}
	with probability at least $1 - \delta$.
\end{theorem}

\begin{proof}[Proof of Theorem~\ref{thm:bias}]
In this proof, when there is no confusion, we use $f(x) = \sum_{i=1}^p e_i(x) f_i$ where $f_i$ denotes the coefficients of $f$ under the basis $e_i(x)$. Adopting this notation, we can write $f(x) = e(x)^* f$ where $f = [f_1, f_2, \ldots, f_p]^T$ also denotes a possibly infinite vector.  	
For the bias, it is easier to work in the frequency domain using the spectral decomposition. Recalling the spectral characterization in the preliminary section,
\begin{align*}
	\bB &= \int \left| e^*(x) T^{1/2} \left[  T^{1/2} e(X) (e(X)^* T e(X))^{-1} e(X)^* T^{1/2} - I \right] T^{-1/2} f_* \right|^2 d\mu(x) \\
	&\leq \int \left\|  \left[  T^{1/2} e(X) (e(X)^* T e(X))^{-1} e(X)^* T^{1/2} - I \right]  T^{1/2} e(x) \right\|^2 d\mu(x) \cdot  \| T^{-1/2} f_* \|^2 \\
	&= \| f_* \|_{\cH}^2  \int \left\|  \left[  T^{1/2} e(X) (e(X)^* T e(X))^{-1} e(X)^* T^{1/2} - I \right]  T^{1/2} e(x) \right\|^2 d\mu(x).
\end{align*}
Here we use the fact that $T^{-1/2} f_* =  \sum_i t_i^{-1/2} f_{*,i} e_i$ and $\| T^{-1/2} f_* \|^2 = \sum_i f_{*,i}^2/t_i = \| f_* \|_{\cH}^2$.
Next, recall the empirical Kernel operator with its spectral decomposition $\widehat{T} = \widehat{U} \widehat{\Lambda} \widehat{U}^*$, with
$\widehat{\Lambda}_{jj} = \frac{1}{n} \lambda_j\left( K(X, X) \right)$. Denote the top $k$ columns of $\widehat{U}$ to be $\widehat{U}_{k}$, and $P_{\widehat{U}_k}^\perp$ to be projection to the eigenspace orthogonal to $\widehat{U}_{k}$. By observing that $T^{1/2} e(X) (e(X)^* T e(X))^{-1} e(X)^* T^{1/2}$ is a projection matrix, it is clear that for all $k\leq n$,
\begin{align}
	\bB &\leq \| f_* \|_{\cH}^2  \int \left\| P^\perp_{\widehat{U}} \left(T^{1/2} e(x) \right) \right\|^2 d\mu(x) \leq  \| f_* \|_{\cH}^2  \int \left\| P^\perp_{\widehat{U}_k} \left(T^{1/2} e(x) \right) \right\|^2 d\mu(x).
\end{align}
We continue the study of the last quantity using techniques inspired by \cite{shawe2004kernel}. Denote the function $g$ indexed by any rank-$k$ projection $U_k$ as
\begin{align}
	g_{U_k}(x) :=  \left\| P_{U_k} \left(T^{1/2} e(x) \right) \right\|^2 = \tr\left(e^*(x) T^{1/2} U_k U_k^T T^{1/2} e(x) \right).
\end{align}
Clearly, $\| U_k U_k^T \|_F = \sqrt{k}$.
Define the function class
\begin{align*}
	\cG_k := \{ g_{U_k}(x): U_k^T U_k = I_k \}.
\end{align*}
It is clear that $g_{\widehat{U}_k} \in \cG_k$. Observe that $g_{\widehat{U}_k}$ is a random function that depends on the data $X$, and we will bound the bias term using the empirical process theory. It is straightforward to verify that
\begin{align*}
	 \E_{\bx \sim \mu} \left\| P^\perp_{\widehat{U}_k} \left(T^{1/2} e(\bx) \right) \right\|^2 &=  \int \left\| P^\perp_{\widehat{U}_k} \left(T^{1/2} e(x) \right) \right\|^2 d\mu(x), \\
	\widehat{\E}_n \left\| P^\perp_{\widehat{U}_k} \left(T^{1/2} e(\bx) \right) \right\|^2 &= \frac{1}{n} \sum_{i=1}^n \left\| P^\perp_{\widehat{U}_k} \left(T^{1/2} e(x_i) \right) \right\|^2 \\
	&=  \tr\left( P^\perp_{\widehat{U}_k} \widehat{T}  P^\perp_{\widehat{U}_k} \right) = \sum_{j > k} \widehat{\Lambda}_{jj} = \frac{1}{n} \sum_{j > k} \lambda_j(K(X, X)).
\end{align*}

Using symmetrization Lemma~\ref{lem:symmetrization} with $M = \sup_{x\in \Omega} K(x, x)$, with probability at least $1-2\delta$,
\begin{align*}
	&  \int \left\| P^\perp_{\widehat{U}_k} \left(T^{1/2} e(x) \right) \right\|^2 d\mu(x) -  \frac{1}{n} \sum_{j > k} \lambda_j(K(X, X))  \\
	= &\E_\mu \left\| P^\perp_{\widehat{U}_k} \left(T^{1/2} e(\bx) \right) \right\|^2  - \widehat{\E}_n \left\| P^\perp_{\widehat{U}_k} \left(T^{1/2} e(\bx) \right) \right\|^2 \\
	\leq & \sup_{U_k: U_k^T U_k = I_k} \left( \E - \widehat{\E}_n \right) \left\| P^\perp_{U_k} \left(T^{1/2} e(\bx) \right) \right\|^2 \\
	\leq & 2\E_\epsilon \sup_{U_k: U_k^T U_k = I_k} \frac{1}{n} \sum_{i=1}^n \epsilon_i \left( \left\| T^{1/2} e(x_i)  \right\|^2 -  \left\| P_{U_k} \left(T^{1/2} e(x_i) \right) \right\|^2 \right) + 3 M
		 \sqrt{\frac{\log 1/\delta}{2n}} 
\end{align*}
by the Pythagorean theorem. Since $\epsilon_i$'s are symmetric and zero-mean and $\left\| T^{1/2} e(x_i)  \right\|^2$ does not depend on $U_k$, the last expression is equal to
\begin{align*}
	& 2\E_\epsilon \sup_{g \in \cG_k} \frac{1}{n} \sum_{i=1}^n \epsilon_i g(x_i)  + 3 M
		 \sqrt{\frac{\log 1/\delta}{2n}}. 
\end{align*}
We further bound the Rademacher complexity of the set $\cG_k$
\begin{align*}
	&\E_\epsilon \sup_{g \in \cG_k} \frac{1}{n}\sum_{i=1}^n \epsilon_i g(x_i) = \E_\epsilon \sup_{U_k} \frac{1}{n}\sum_{i=1}^n \epsilon_i g_{U_k}(x_i) \\
	& = \E_\epsilon \frac{1}{n} \sup_{U_k} \left\langle U_k U_k^T,  \sum_{i=1}^n \epsilon_i T^{1/2} e(x_i) e^*(x_i) T^{1/2}  \right\rangle \\
	& \leq \frac{\sqrt{k}}{n} \E_\epsilon \left\| \sum_{i=1}^n \epsilon_i T^{1/2} e(x_i) e^*(x_i) T^{1/2}  \right\|_F 
\end{align*}
by the Cauchy-Schwarz inequality and the fact that $\| U_k U_k^T\|_F \leq \sqrt{k}$. The last expression is can be further evaluated by the independence of $\epsilon_i$'s
\begin{align*}	
	\frac{\sqrt{k}}{n} \left\{ \E_\epsilon \left\| \sum_{i=1}^n \epsilon_i T^{1/2} e(x_i) e^*(x_i) T^{1/2}  \right\|_F^2 \right\}^{1/2} & = \frac{\sqrt{k}}{n} \left\{ \sum_{i=1}^n \left\|  T^{1/2} e(x_i) e^*(x_i) T^{1/2}  \right\|_F^2  \right\}^{1/2} \\
	& = \sqrt{\frac{k}{n}} \sqrt{ \frac{\sum_{i=1}^n K(x_i, x_i)^2}{n}}.
\end{align*}

Therefore, for all $k\leq n$, with probability at least $1 - 2n\delta$,
\begin{align*}
	\bB \leq \| f_* \|_{\cH}^2 \cdot \inf_{0\leq k \leq n}  \left\{  \frac{1}{n} \sum_{j > k} \lambda_j(K(X, X)) + 2 \sqrt{\frac{k}{n}} \sqrt{ \frac{\sum_{i=1}^n K(x_i, x_i)^2}{n}} + 3M \sqrt{\frac{\log 1/\delta}{2n}} \right\}.
\end{align*}

\end{proof}

\begin{remark}
	Let us compare the bounds obtained in this paper to those one can obtain for classification with a margin. For classification, Thm. 21 in
	\cite{bartlett2002rademacher} shows that the misclassification error is upper bounded with probability at least $1 - \delta$ as 
	\begin{align*}
		\E {\bf 1}(\by \widehat{f}(\bx)<0) \leq \E \phi_\gamma(\by \widehat{f}(\bx)) \leq \widehat{\E}_n \phi_\gamma(\by \widehat{f}(\bx)) + \frac{C_\delta}{\gamma \sqrt{n}} \sqrt{\frac{\sum_{i=1}^n K(x_i, x_i)}{n}}
	\end{align*}
	where $\phi_{\gamma}(t) := \max(0, 1-t/\gamma) \wedge 1$ is the margin loss surrogate for the indicator loss ${\bf 1}(t<0)$. By tuning the margin $\gamma$, one obtains a family of upper bounds. 
	
	Now consider the noiseless regression scenario (i.e. $\sigma=0$ in (A.1)). In this case, the variance contribution to the risk is zero, and
	\begin{align*}
		\E_{Y|X} \| \widehat{f} -\by \|_{L^2_\mu}^2 &= \E_{Y|X} \| \widehat{f} - f_*\|^2_{L^2_\mu} = \E [ P_n^{\perp} f_* ]^2 \leq \E [P_k^{\perp}f_*]^2   \\
		& \leq \widehat{\E}_n [P_k^{\perp}f_*]^2 + C'_\delta \sqrt{\frac{k}{n}} \sqrt{ \frac{\sum_{i=1}^n K(x_i, x_i)^2}{n}}
	\end{align*}
	where $P_k$ is the best-rank $k$ projection (based on $X$) and $P_k^\perp$ denotes its orthogonal projection. By tuning the parameter $k$ (similar as the $1/\gamma$ in classification), one can balance the RHS to obtain the optimal trade-off. 
	
	However, classification is easier than regression in the following sense: $\widehat{f}$ can present a non-vanishing bias in estimating $f_*$, but as long as the bias is below the empirical margin level, it plays no effect in the margin loss $\phi_{\gamma}(\cdot)$. In fact, for classification, under certain conditions, one can prove exponential convergence for the generalization error \citep{koltchinskii2005exponential}.
\end{remark}

\section{Experiments}
\label{sec:experiments}

\subsection{MNIST}

In this section we provide full details of the experiments on MNIST \citep{lecun2010mnist}. Our first experiment considers the following problem: for each pair of distinct digits $(i, j)$, $i, j \in \{0,1,\ldots,9\}$, label one digit as $1$ and the other as $-1$, then fit the Kernel Ridge Regression with Gaussian kernel $k(x, x') = \exp(-\|x - x' \|^2/d)$, where $d=784$ is the dimension as analyzed in our theory (also the default choice in Scikit-learn package \citep{scikit-learn}). For each of the $\binom{10}{2} = 45$ pairs of experiments, we chose $\lambda = 0$ (no regularization, interpolation estimator), $\lambda = 0.1$ and $\lambda = 1$. We evaluated the performance on the \emph{out-of-sample} test dataset, with the error metric
\begin{align}
	\label{eq:normalized_mse}
	\frac{\sum_{i} (\widehat{f}(x_i) - y_i)^2}{\sum_{i} (\bar{y} - y_i)^2}.
\end{align}
Remarkably, among all 45 experiments, no-regularization performs the best. We refer to the table in Section~\ref{sec:mnist} for a complete list of numerical results. For each experiment, the sample size is roughly $n\approx 10000$.

\begin{figure}[h]
\centering
\includegraphics[width=0.8\textwidth]{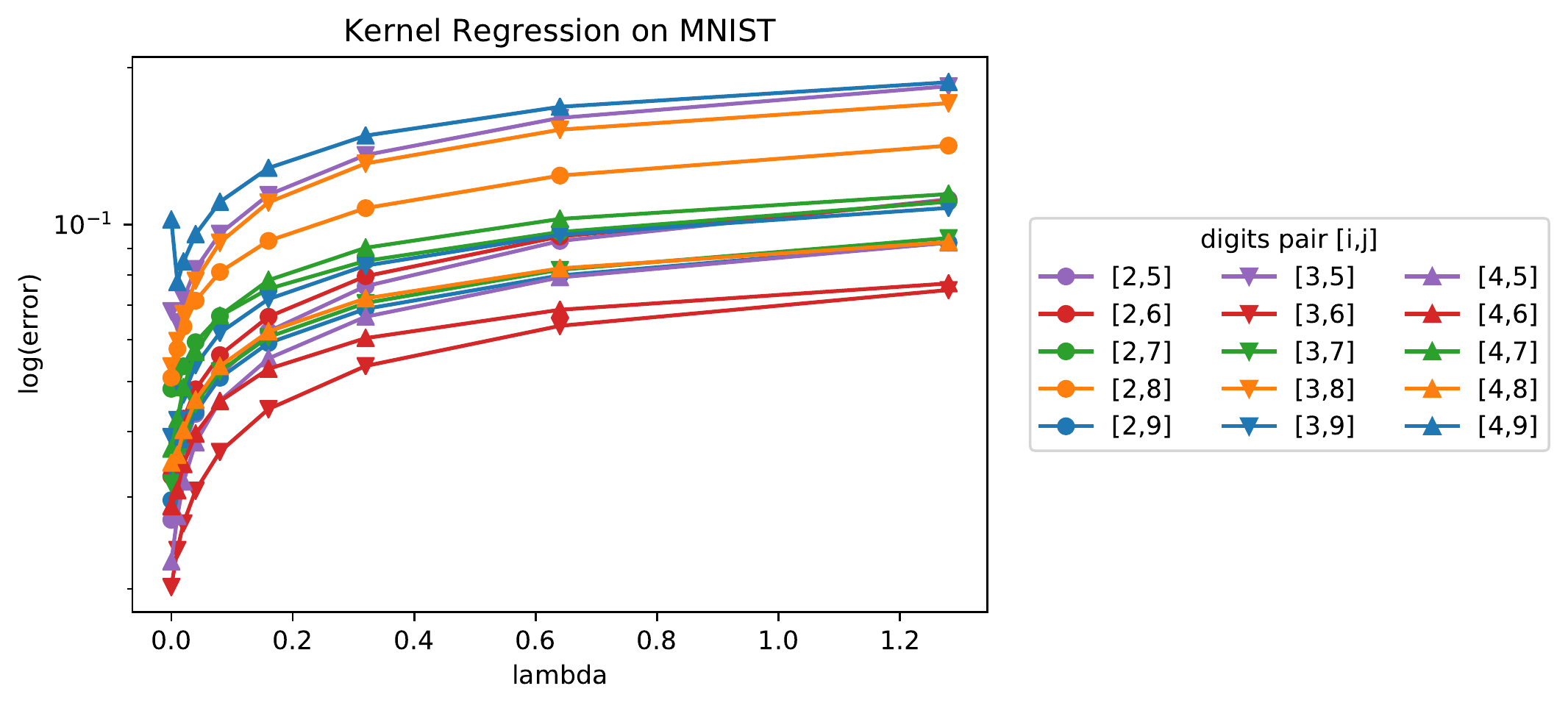}
\caption{Test error, normalized as in \eqref{eq:normalized_mse}. The y-axis is on the log scale.}
\label{fig:mnisit-finer-grid}
\end{figure}

The second experiment is to perform the similar task on a finer grid of regularization parameter $\lambda \in \{0, 0.01, 0.02, 0.04, 0.08, 0.16, 0.32, 0.64, 1.28\}$. Again, in all but one pair, the interpolation estimator performs the best in out-of-sample prediction. We refer to Figure~\ref{fig:mnisit-finer-grid} for details.


\begin{figure}[ht]
\centering
\includegraphics[width=0.45\textwidth]{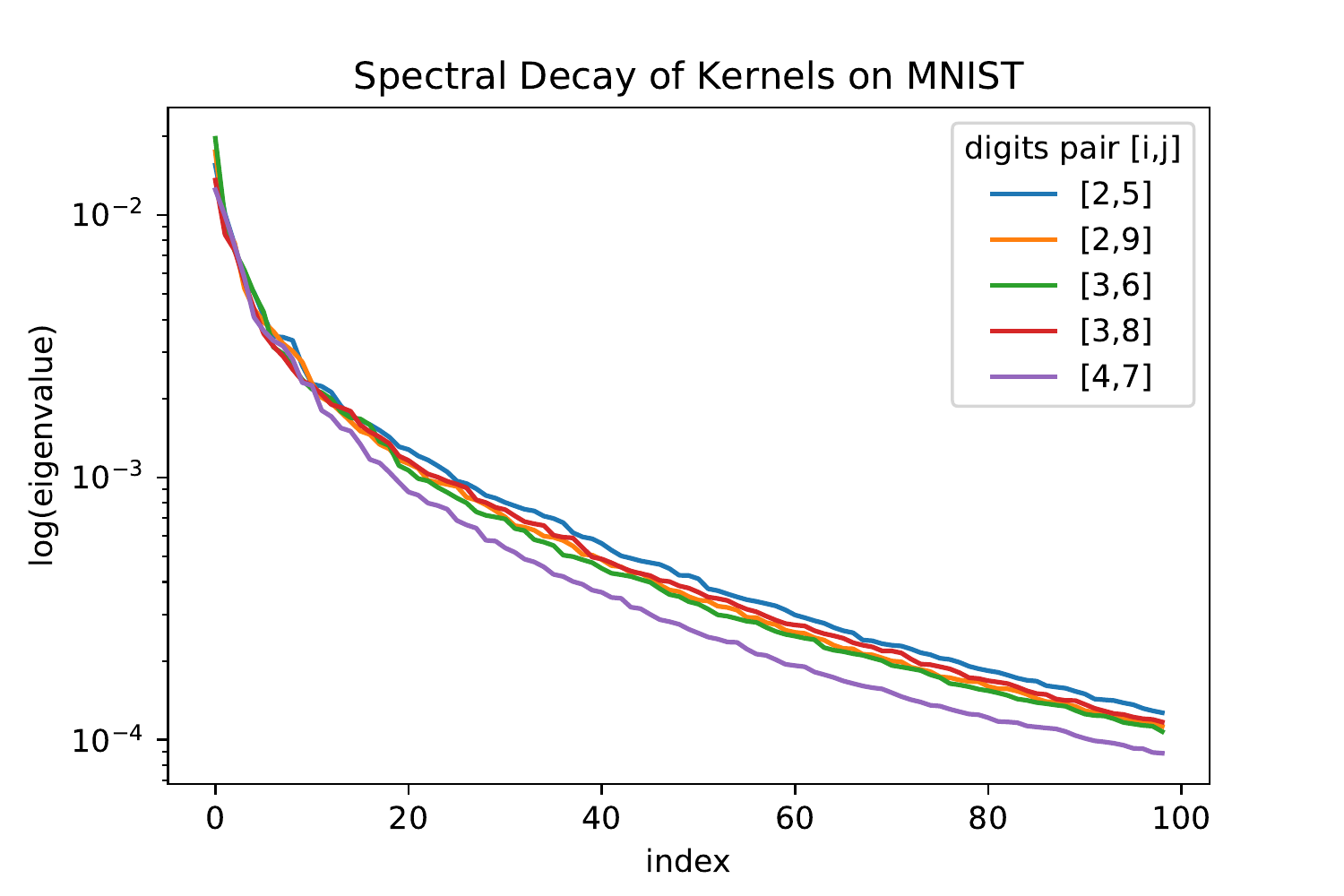}
\includegraphics[width=0.45\textwidth]{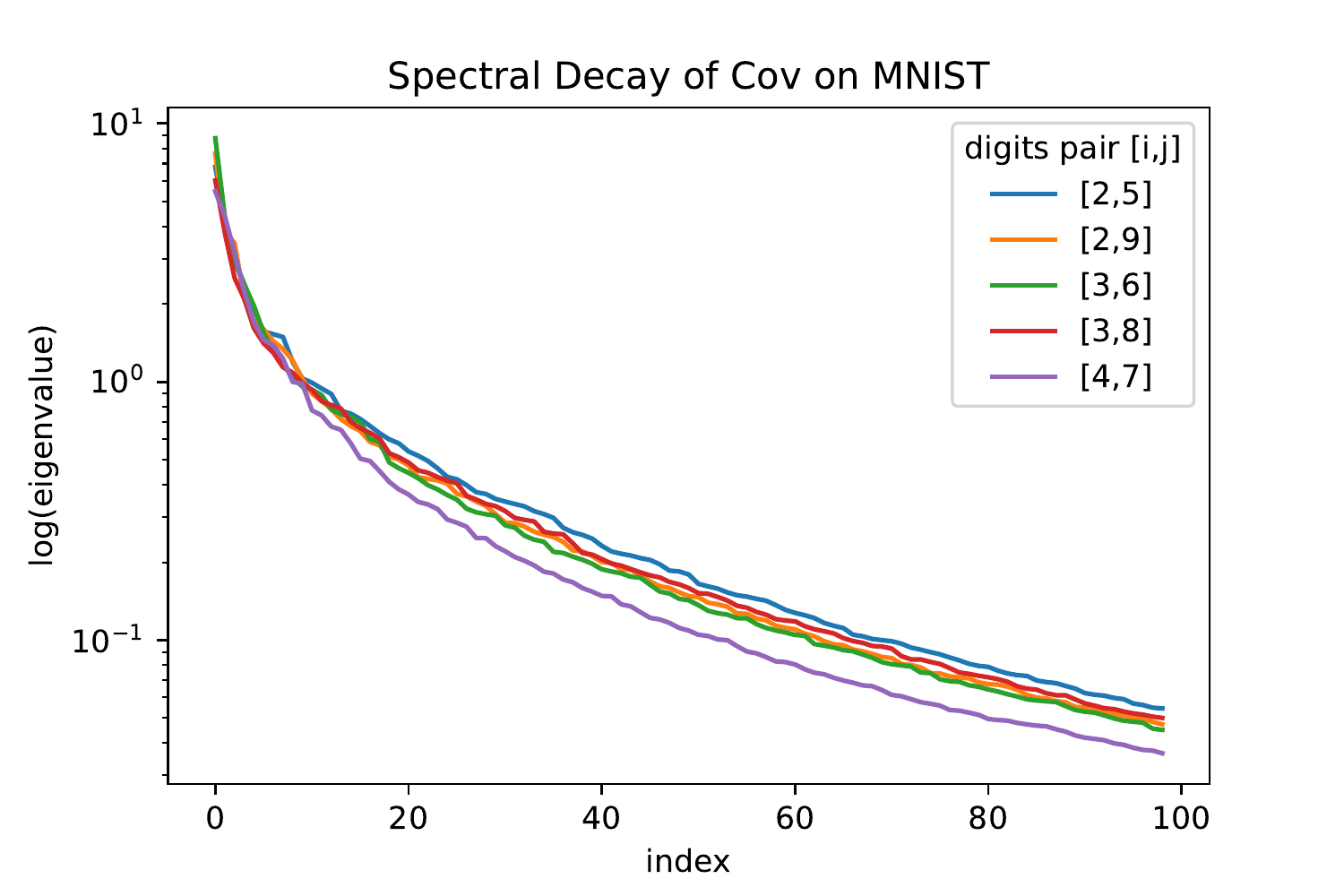}
\caption{Spectral decay. The y-axis is on the log scale.}
\label{fig:mnisit-spectral}
\end{figure}

To conclude this experiment, we plot the eigenvalue decay of the empirical kernel matrix and the sample covariance matrix for the 5 experiments shown in the introduction. The two plots are shown in Figure~\ref{fig:mnisit-spectral}. Both plots exhibit a fast decay of eigenvalues, supporting the theoretical finding that interpolation performs well on a test set in such situations. 

On the other hand, it is easy to construct examples where the eigenvalues do not decay and interpolation performs poorly. This is the case, for instance, if $X_i$ are i.i.d. from spherical Gaussian. One can show that in the high-dimensional regime, the variance term itself (and not just the upper bound on it) is large. Since the bias-variance decomposition is an equality, it is not possible to establish good $L^2_\mu$ convergence.

\subsection{A Synthetic Example}
\label{sec:synthetic}

\begin{figure}[ht]
\centering
\includegraphics[width=0.32\textwidth]{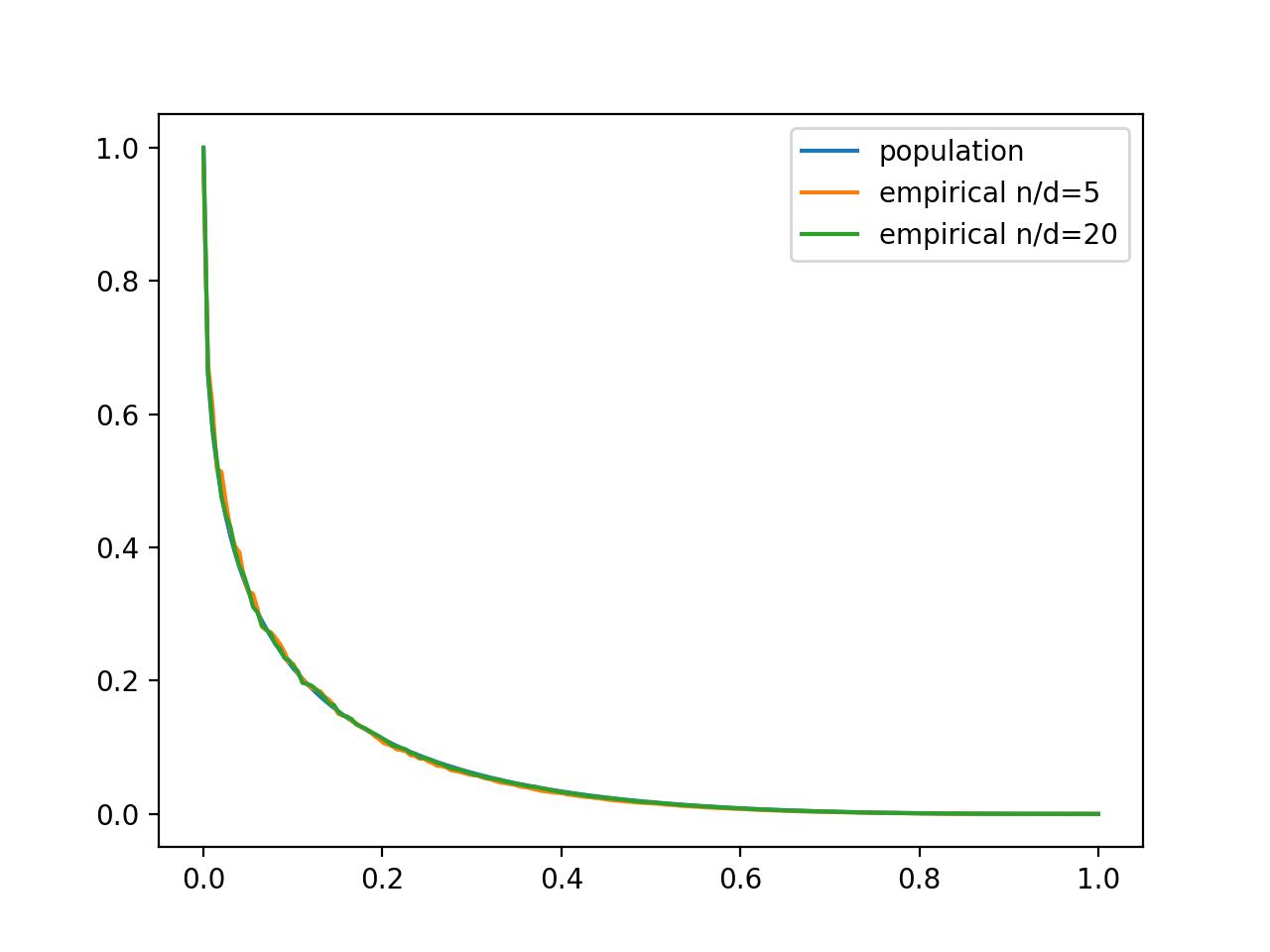}
\includegraphics[width=0.32\textwidth]{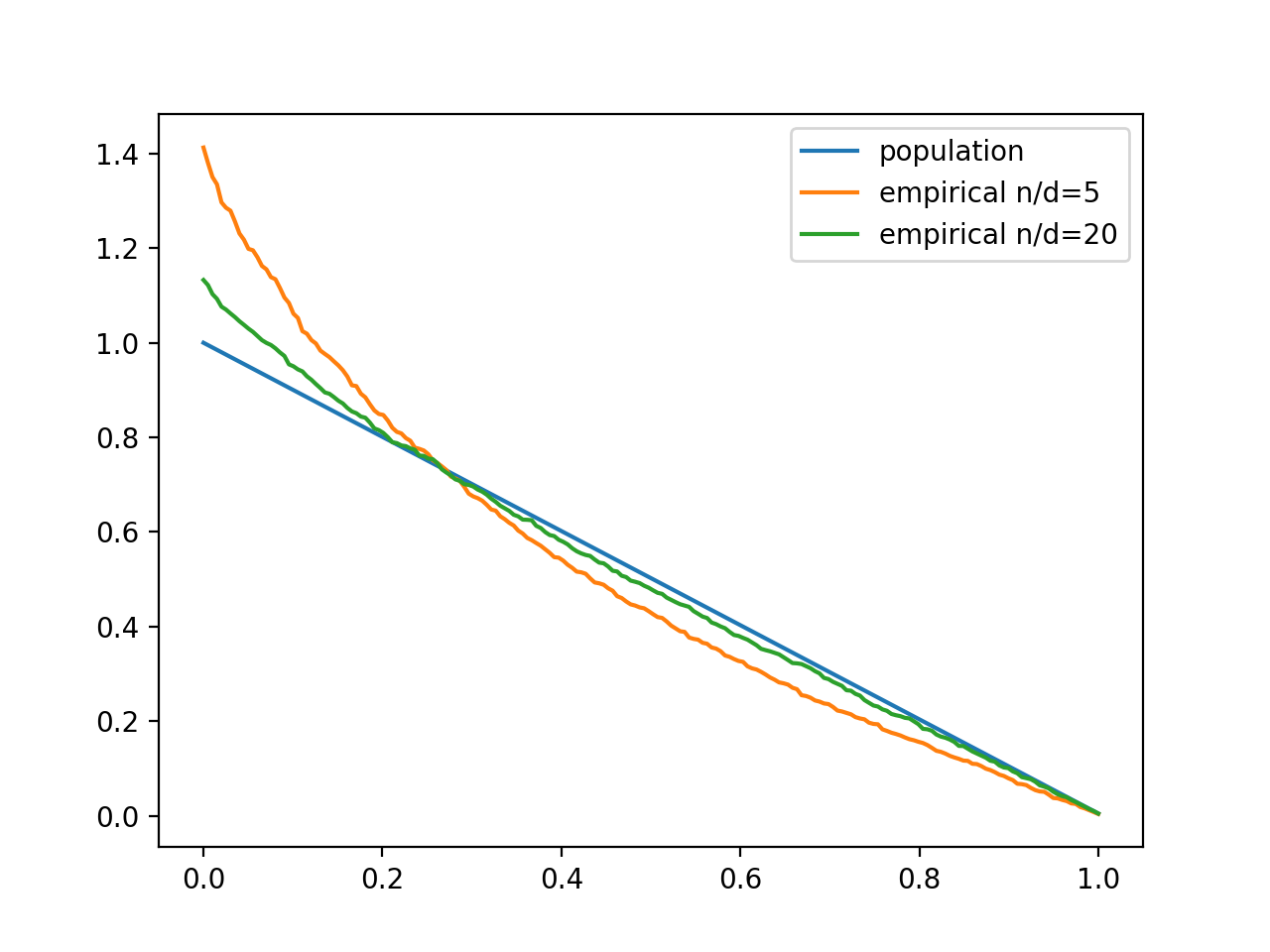}
\includegraphics[width=0.32\textwidth]{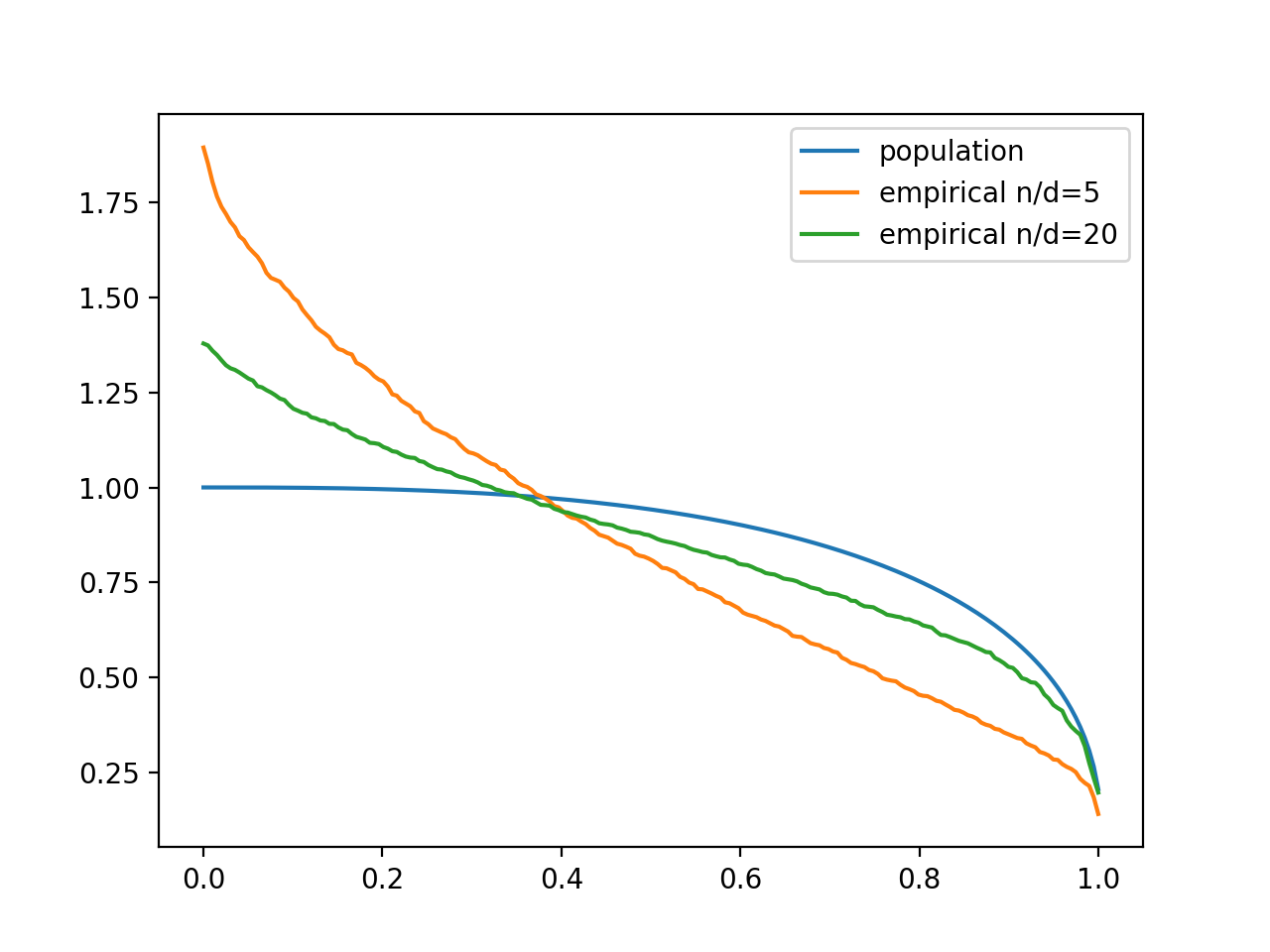}

\includegraphics[width=0.32\textwidth]{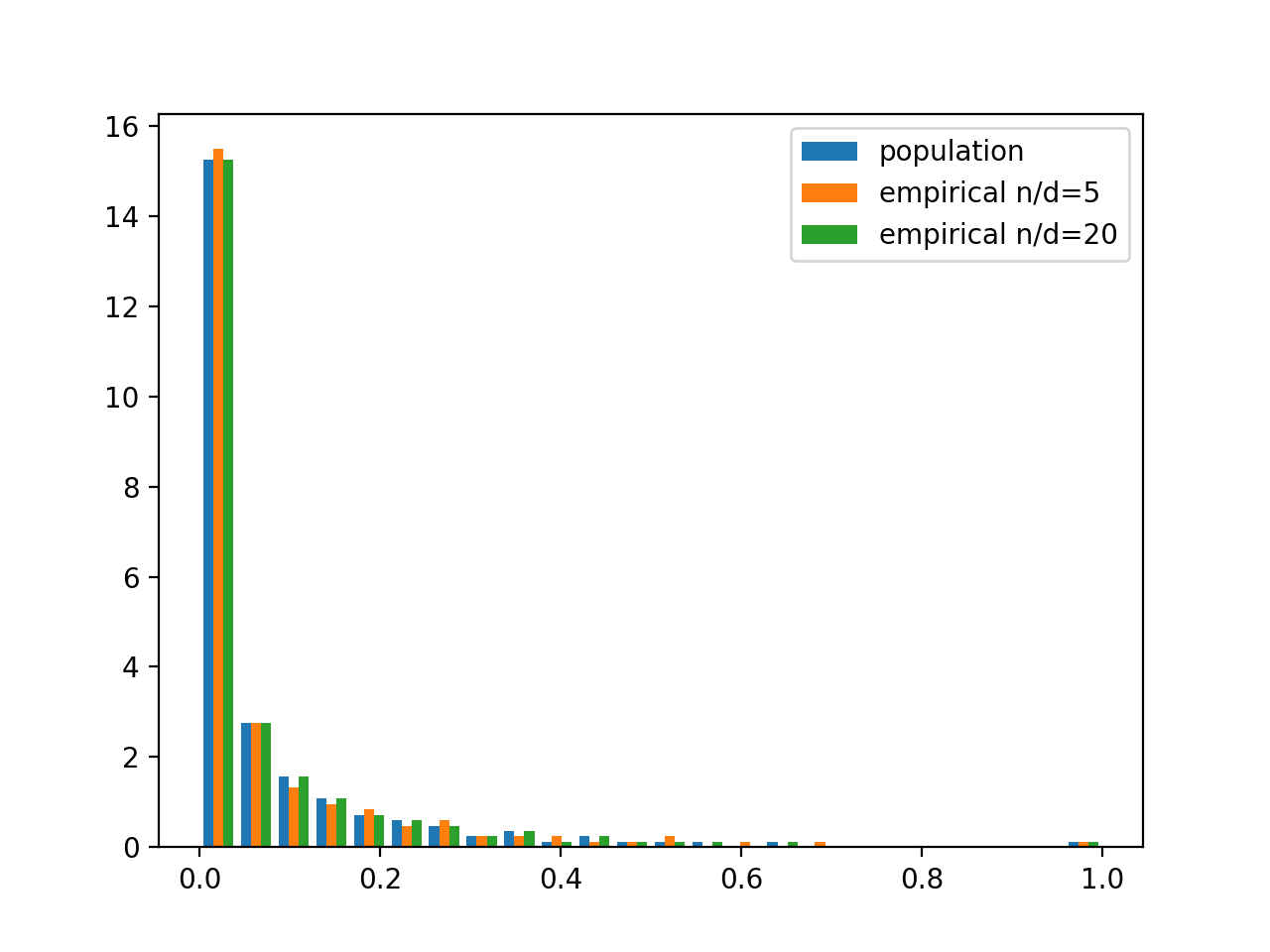}
\includegraphics[width=0.32\textwidth]{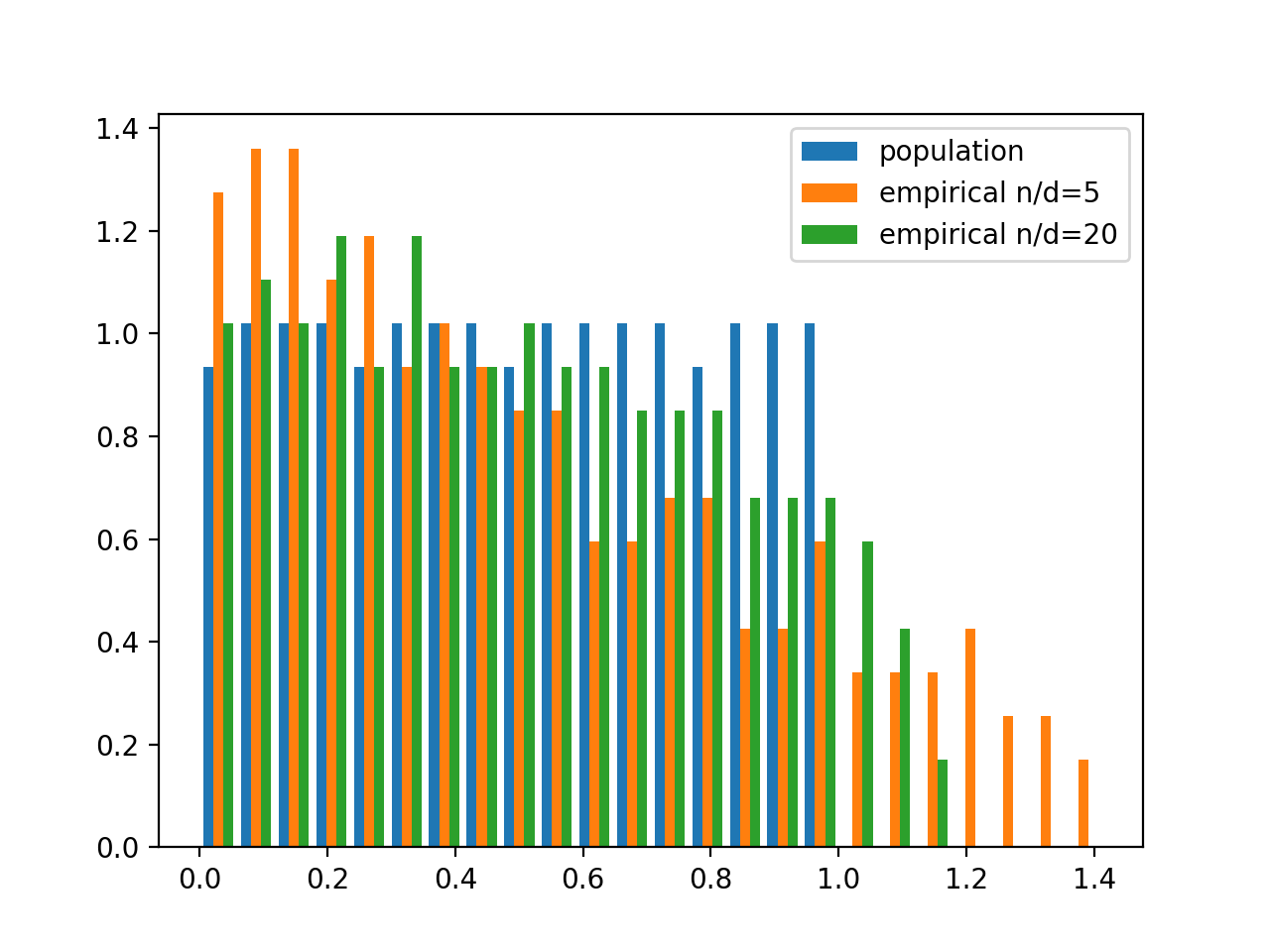}
\includegraphics[width=0.32\textwidth]{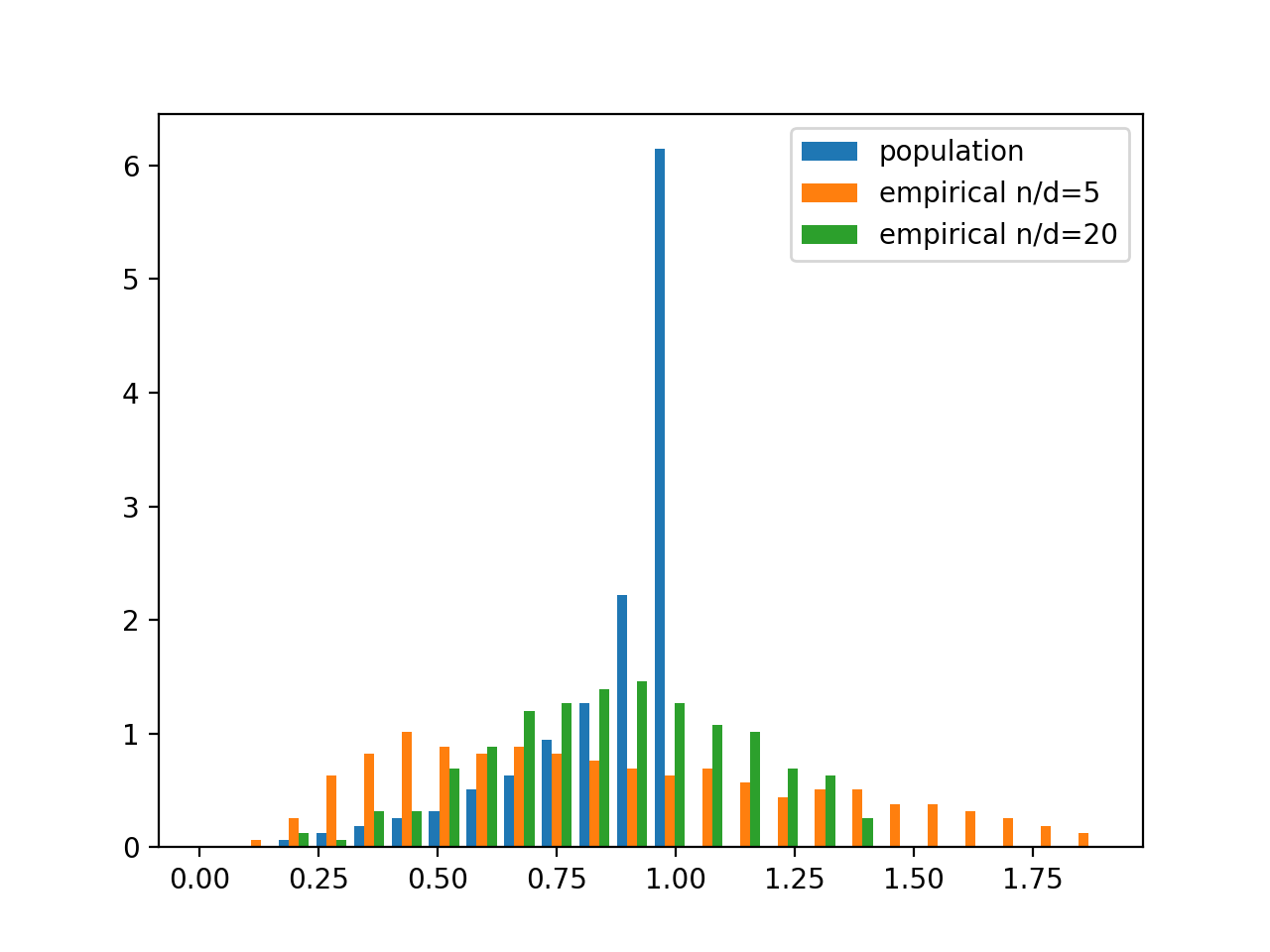}

\caption{Varying spectral decay: case $n > d$. Columns from left to right: $\kappa = e^{-1}, e^{0}, e^{1}$. Rows from top to bottom: ordered eigenvalues, and the histogram of eigenvalues. Here we plot the population eigenvalues for $\Sigma_d$, and the empirical eigenvalues for $X^*X/n$. In this simulation, $d = 100$, $n = 500, 2000$. }
\label{fig:n_m_d}
\end{figure}

\begin{figure}[ht]
\centering
\includegraphics[width=0.32\textwidth]{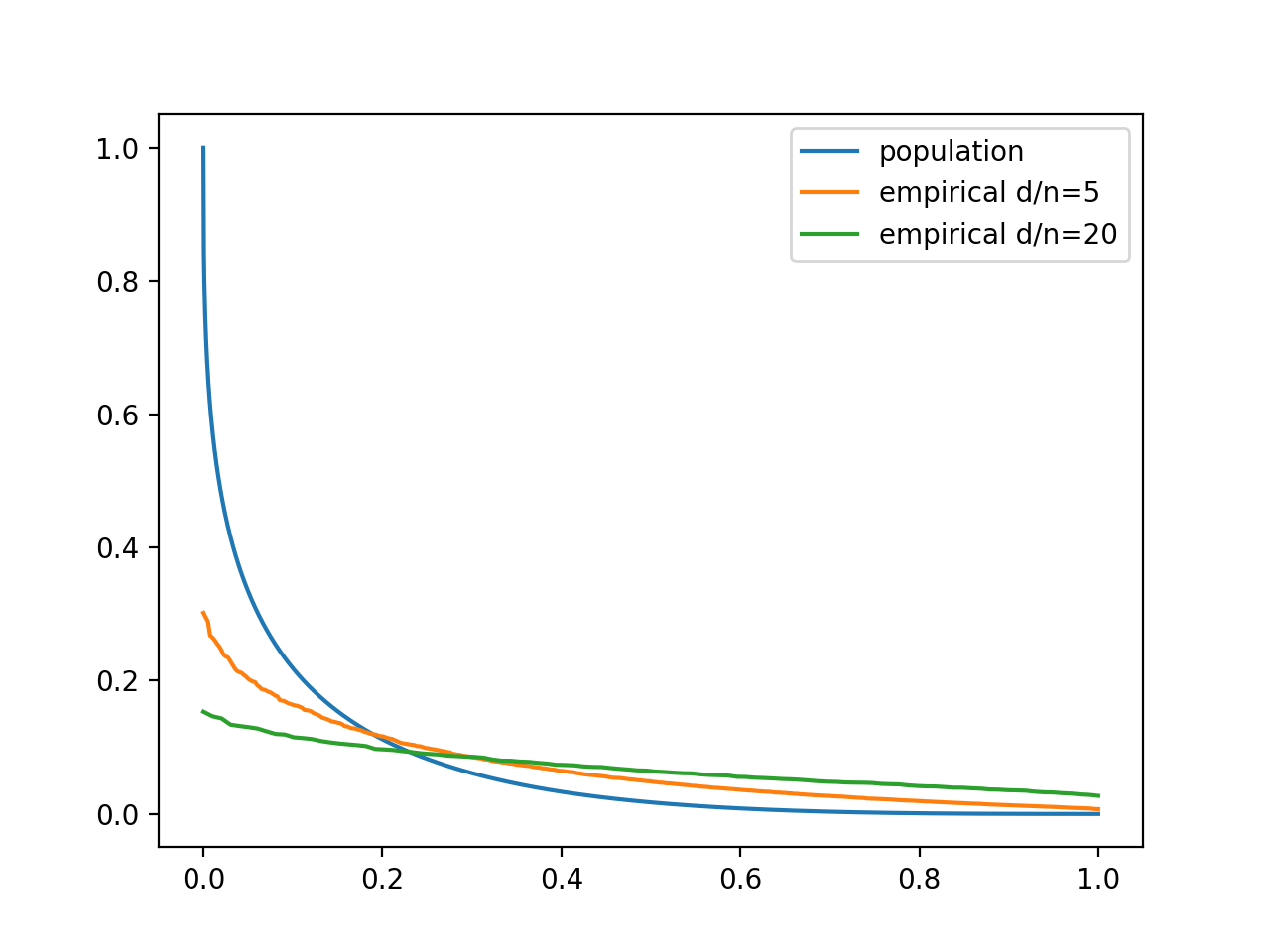}
\includegraphics[width=0.32\textwidth]{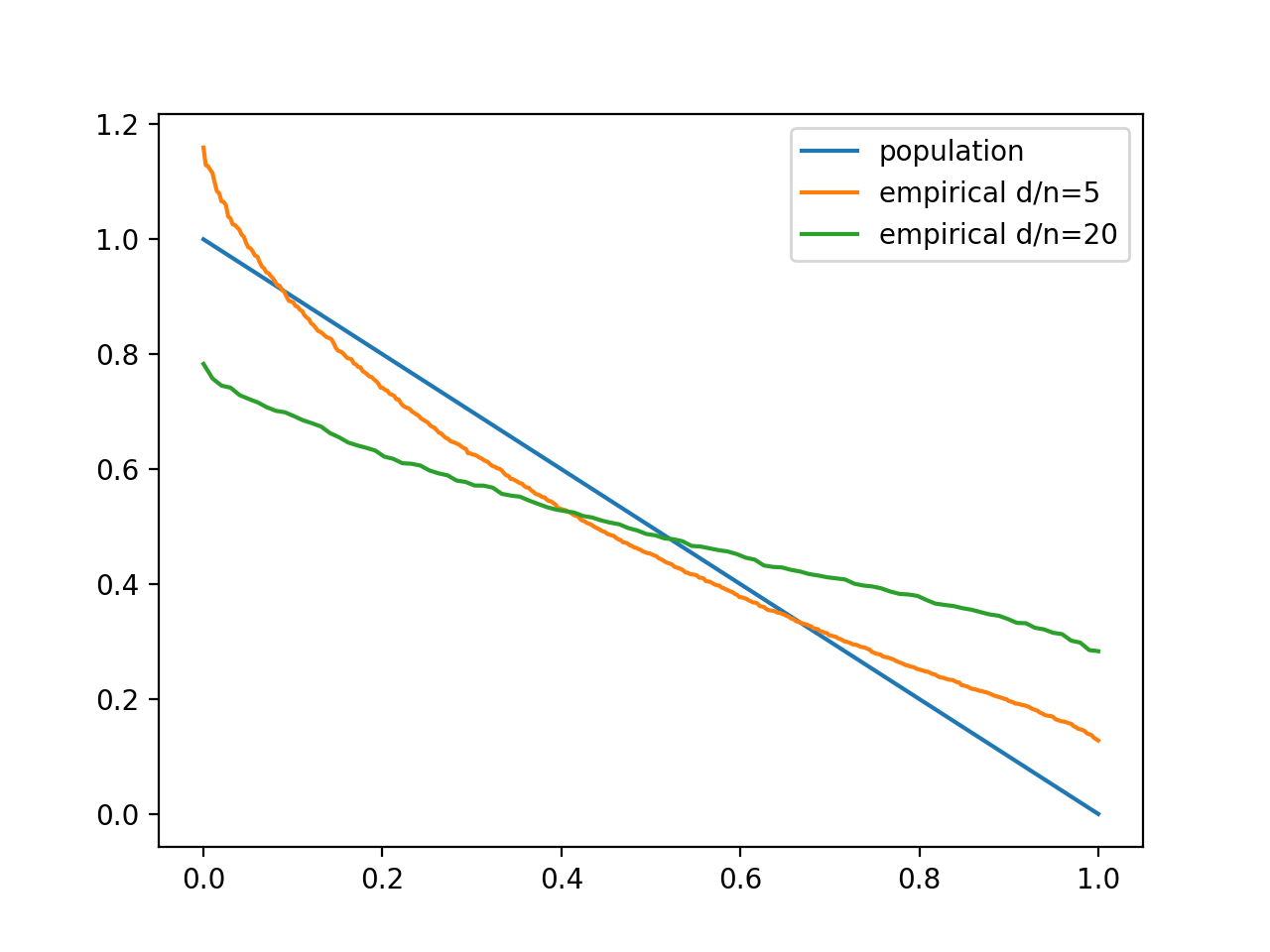}
\includegraphics[width=0.32\textwidth]{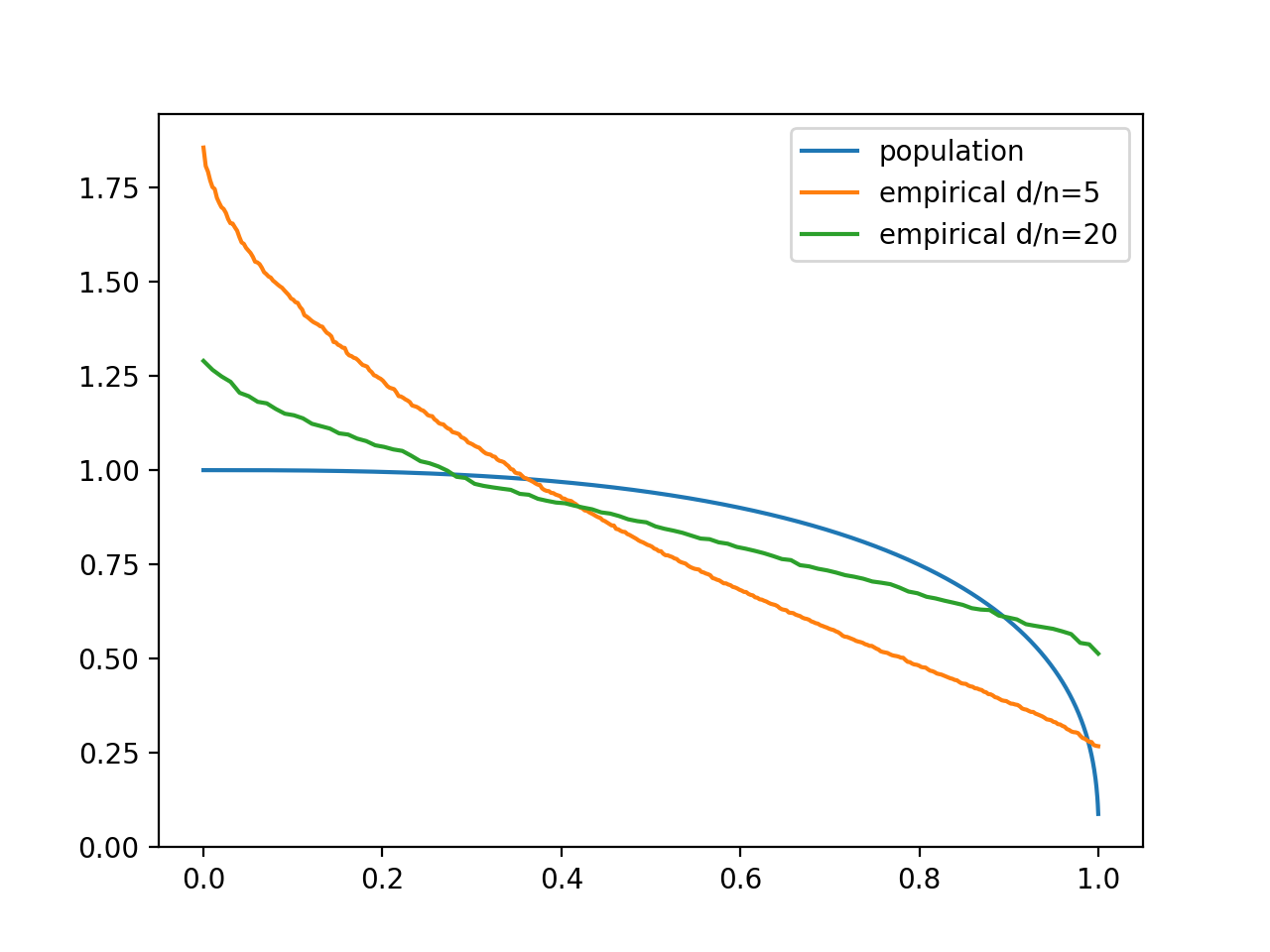}

\includegraphics[width=0.32\textwidth]{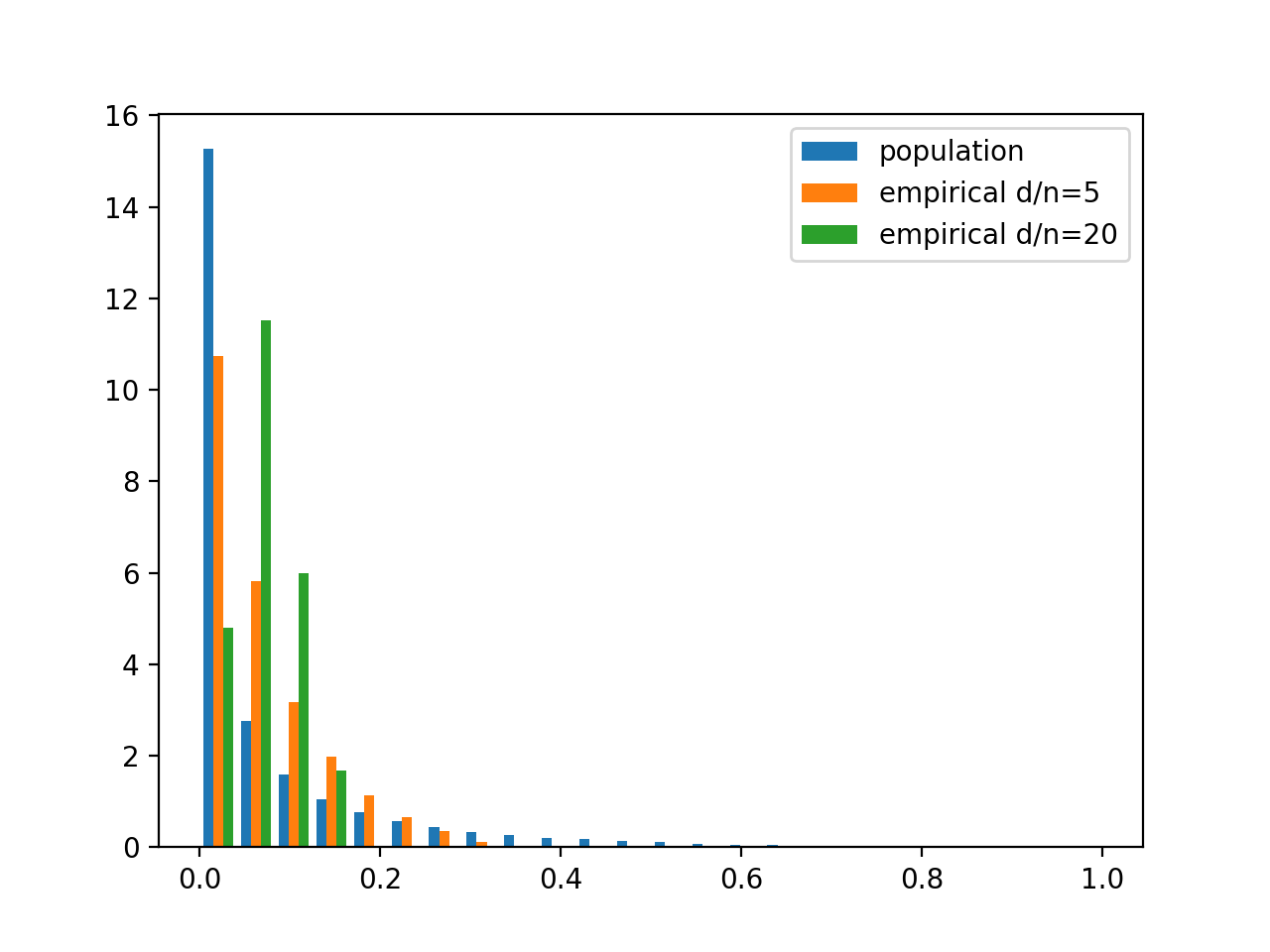}
\includegraphics[width=0.32\textwidth]{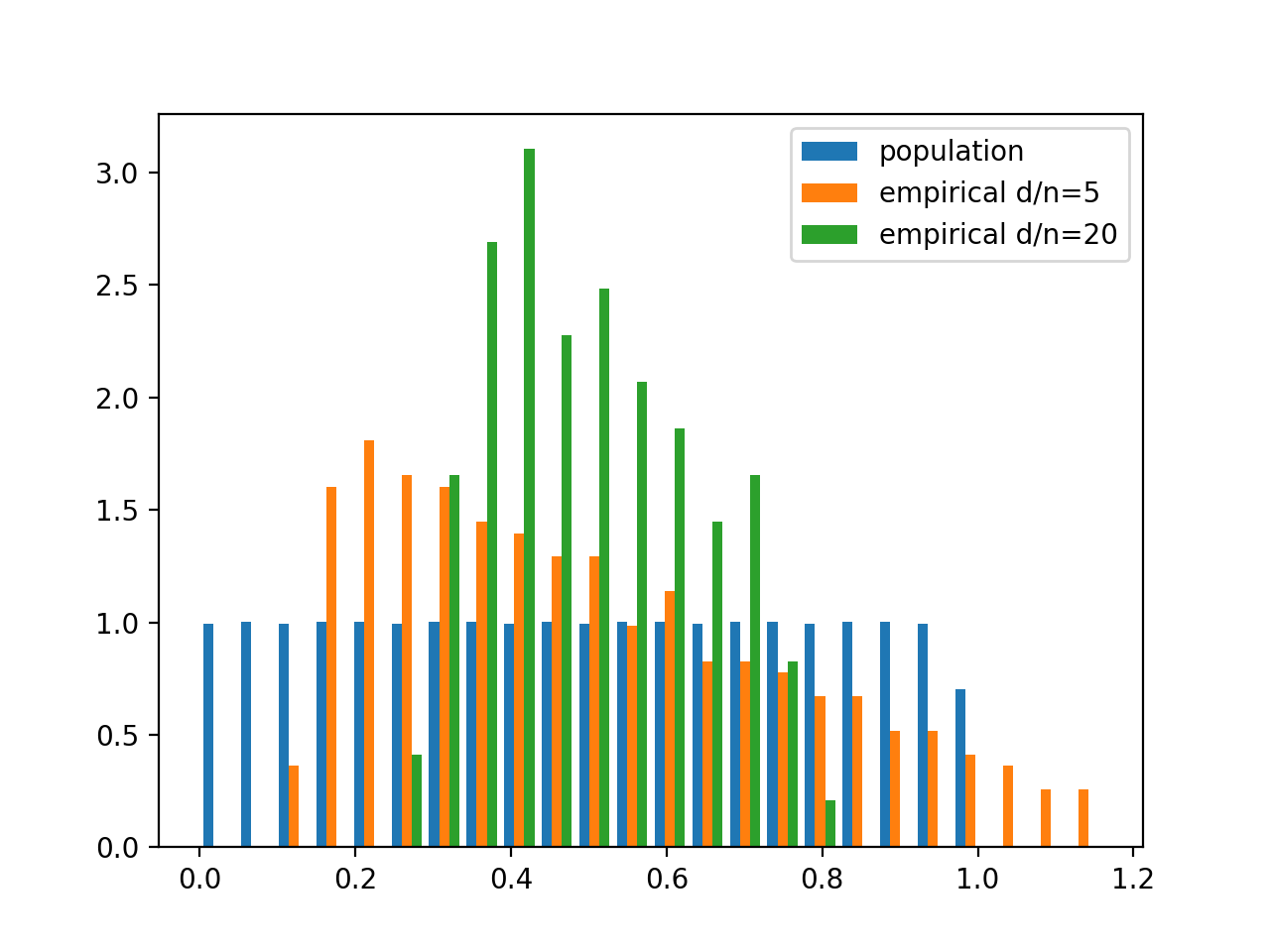}
\includegraphics[width=0.32\textwidth]{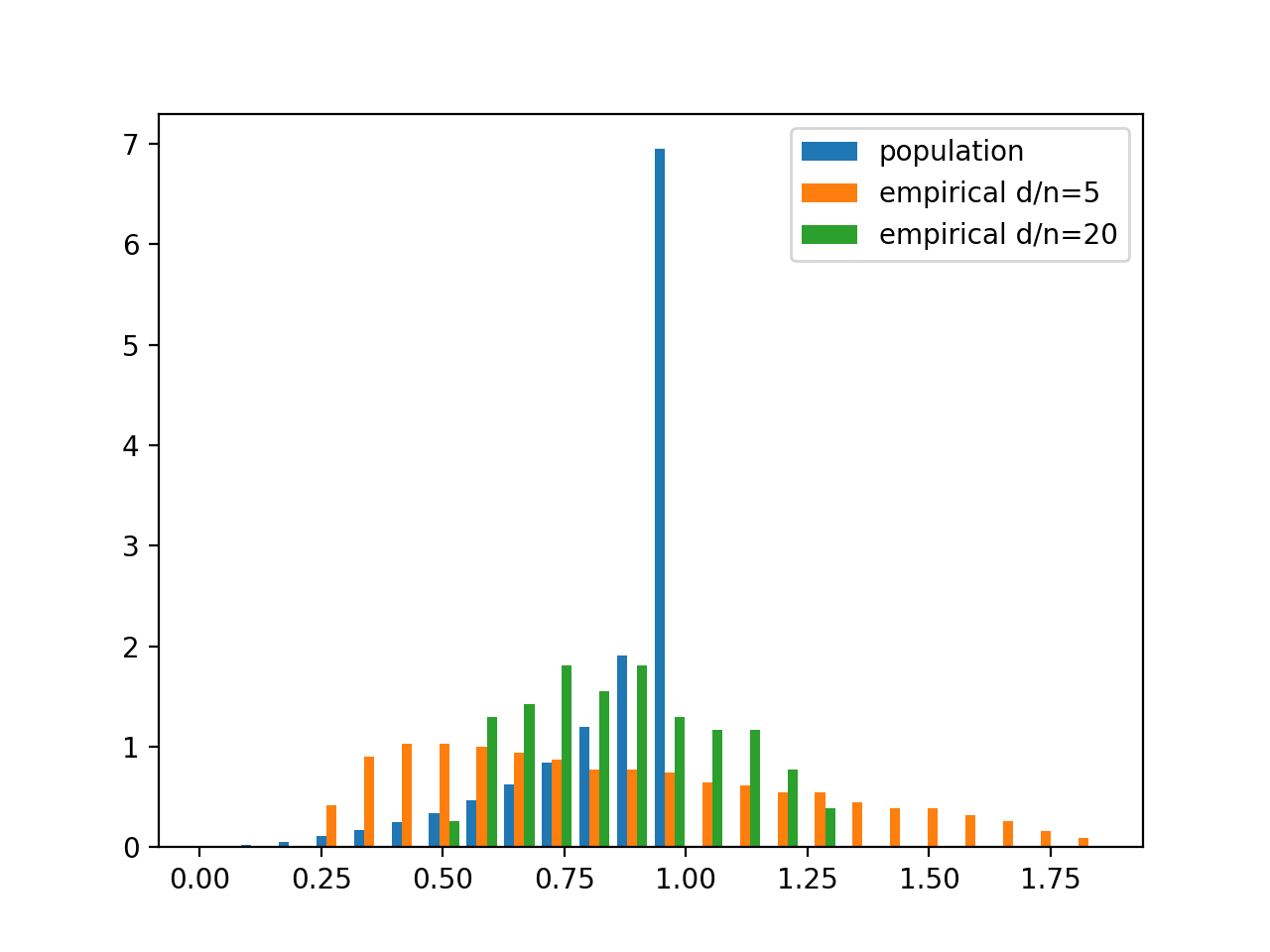}

\caption{Varying spectral decay: case $d > n$. Columns from left to right: $\kappa = e^{-1}, e^{0}, e^{1}$. Rows from top to bottom: ordered eigenvalues, and the histogram of eigenvalues. Here we plot the population eigenvalues for $\Sigma_d$, and the empirical eigenvalues for $XX^*/d$. In this simulation, $d = 2000$, $n = 400, 100$. }
\label{fig:d_m_n}
\end{figure}

\begin{figure}[H]
			\centering
			\includegraphics[width=0.45\textwidth]{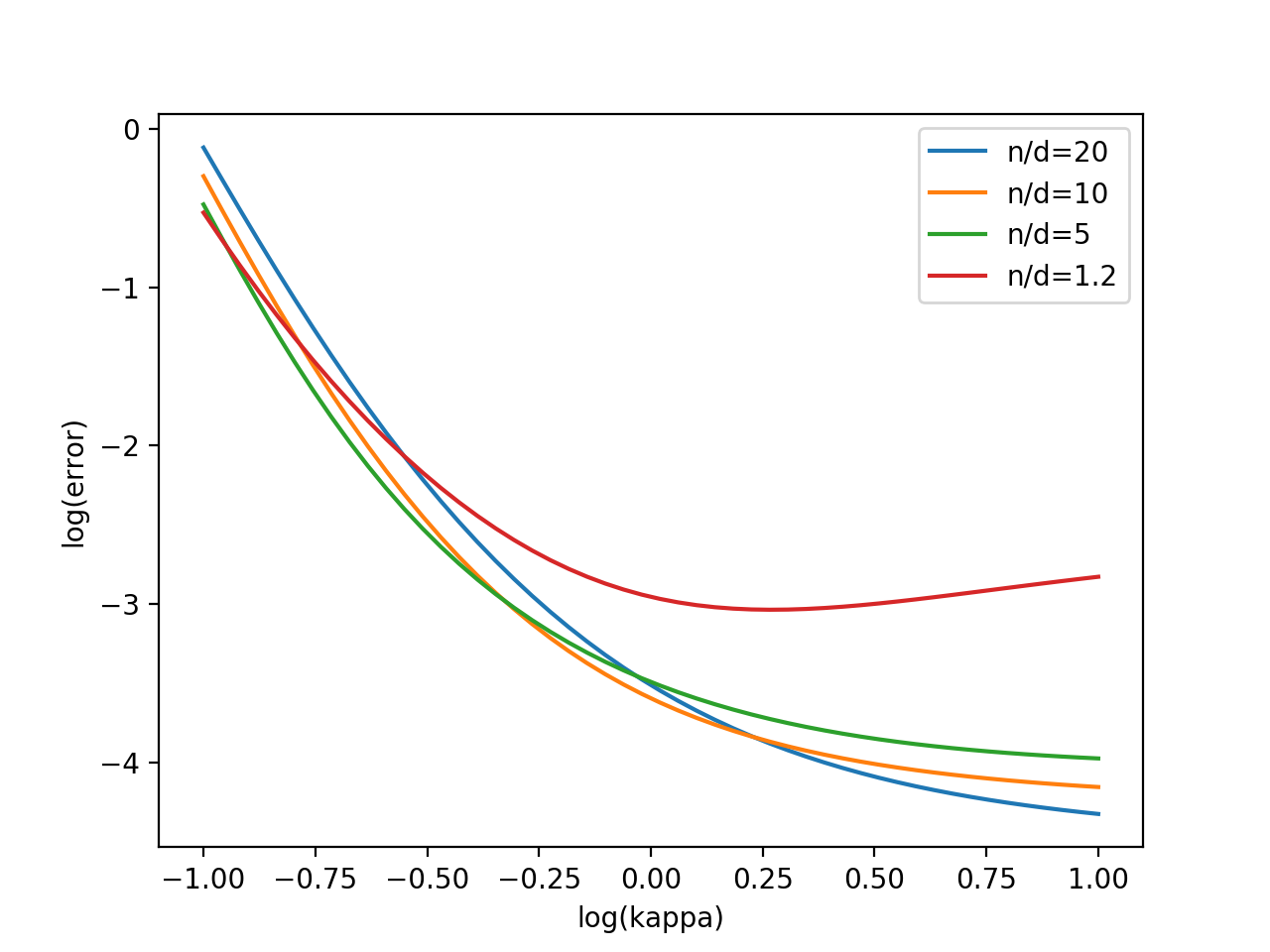}
			\includegraphics[width=0.45\textwidth]{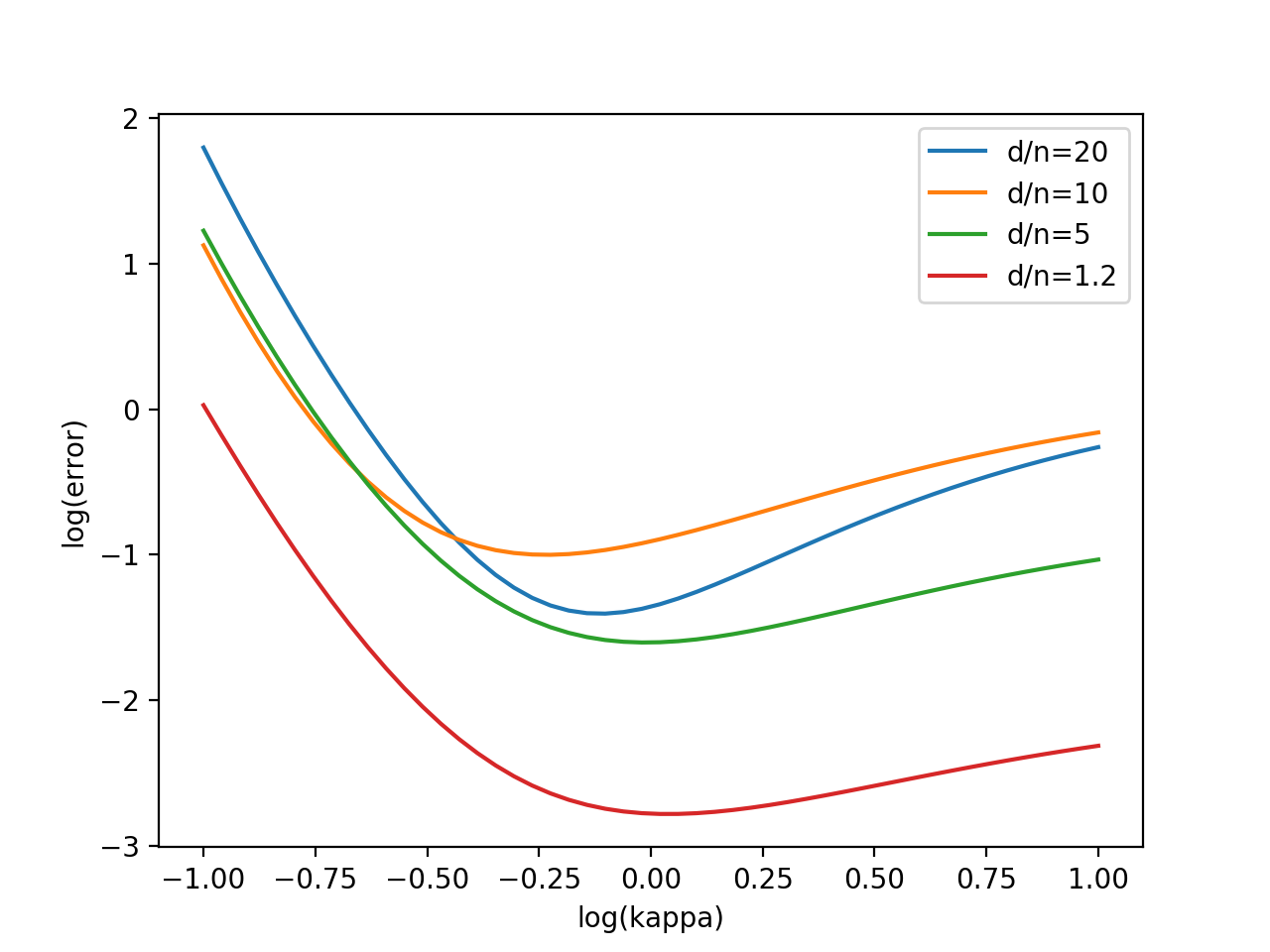}

			\caption{Varying spectral decay: generalization error for high noise case. Left: $d = 200$, $n = 4000, 2000, 1000, 240$. Right: $n = 200$, $d = 4000, 2000, 1000, 240$.}
			\label{fig:gen_error_high_noise}
			\end{figure}

In this section we provide the details of the synthetic experiments mentioned in Section~\ref{sec:data-dependent-bound} for Tables~\ref{tab:simuls_n_d}-\ref{tab:simuls_d_n} and Figures~\ref{fig:gen_error_n_d}-\ref{fig:gen_error_d_n}. We choose the RBF kernel as the non-linearity with $h(t)=\exp(-t)$. Again, we consider a family of eigenvalue decays for the covariance matrix parametrized by $\kappa$, with the small $\kappa$ describing fast spectral decay
\begin{align*}
	\lambda_j(\Sigma_{d, \kappa}) = \left(1 - ((j-1)/d)^{\kappa}\right)^{1/\kappa},  1\leq j \leq d.
\end{align*}
We set a target non-linear function $f_*$ in the RKHS with kernel $K(x, x') = h\left(\| x - x' \|^2/d \right)$ as
\begin{align*}
	f_*(x) = \sum_{l=1}^{100} K(x, \theta_l), ~~\theta_l \stackrel{i.i.d.}{\sim} N(0, I_d).
\end{align*}
For each parameter triplet $(n, d, \kappa)$, we generate data in the following way
\begin{align*}
	x_i \sim N(0, \Sigma_{d, \kappa}), ~~ y_i = f_*(x_i) + \epsilon_i
\end{align*}
for $1\leq i \leq n$ where $\epsilon_i \sim N(0, \sigma^2)$ is independent noise, with $\sigma = 0.1$ (Figures~\ref{fig:gen_error_n_d}-\ref{fig:gen_error_d_n}) and $\sigma = 0.5$ (Figures~\ref{fig:gen_error_high_noise}). 
 Figures~\ref{fig:n_m_d}-\ref{fig:d_m_n} contrasts the difference between the population and empirical eigenvalues for various parameter triplets $(n, d, \kappa)$.

We now explain Figures~\ref{fig:gen_error_n_d}-\ref{fig:gen_error_d_n}, which illustrate the true generalization error in this synthetic example, by varying the spectral decay $\kappa$, for a particular case of high dimensionality ratio $d/n$. Here we plot the \textit{out-of-sample} test error for the interpolated min-norm estimator $\widehat{f}$ on fresh new test data $(x_t, y_t)$ from the same data generating process, with the error metric
\begin{align*}
	\text{error} = \frac{\sum_t (\widehat{f}(x_t) - f_*(x_t) )}{\sum_t ( y_t - \bar{y} )^2}.
\end{align*}
The error plots are shown in Figure~\ref{fig:gen_error_n_d} (for $n>d$) and \ref{fig:gen_error_d_n} (for $d>n$), and Figure~\ref{fig:gen_error_high_noise} for the high noise case. On the x-axis, we plot the $\log(\kappa)$, and on the y-axis the $\log(\text{error})$. Each curve corresponds to the generalization error behavior (and the bias and variance trade-off) as we vary spectral decay from fast to slow (as $\kappa$ increases) for a particular choice of $d/n$ or $n/d$ ratio. Clearly, for a general pair of high dimensionality ratio $d/n$, there is a ``sweet spot'' of $\kappa$ (favorable geometric structure) such that the trade-off is optimized.

\section{Further Discussion}

This paper is motivated by the work of \cite{belkin2018understand} and  \cite{zhang2016understanding}, who, among others, observed the good out-of-sample performance of interpolating rules. This paper continues the line of work in \citep{belkin2018overfitting,belkin2018does,belkin2018approximation} on understanding theoretical mechanisms for the good out-of-sample performance of interpolation. We leave further investigations on the connection between kernel ridgeless regression and two-layer neural networks as a future work \citep{dou2019training}.

From an algorithmic point of view, the minimum-norm interpolating solution can be found either by inverting the kernel matrix, or by performing gradient descent on the least-squares objective (starting from $0$). Our analysis can then be viewed in the light of recent work on implicit regularization of optimization procedures \citep{yao2007early, neyshabur2014search,gunasekar2017implicit,li2017algorithmic}.

The paper also highlights a novel type of implicit regularization. In addition, we discover that once we parametrize the geometric properties --- the spectral decay --- we discover the familiar picture of the bias-variance trade-off, controlled by the implicit regularization that adapts to the favorable geometric property of the data.
Moreover, if one explicitly parametrizes the choice of the kernel by, say, the bandwidth, we are likely to see the familiar picture of the bias-variance trade-off, despite the fact that the estimator is always interpolating. 
Whether one can achieve optimal rates of estimation (under appropriate assumptions) for the right choice of the bandwidth appears to be an interesting and difficult statistical question. Another open question is whether one can characterize situations when the interpolating minimum-norm solution is dominating the regularized solution in terms of expected performance.


\bibliography{ref}

\newpage

\appendix

\section{Propositions}

We first borrow a technical result for concentration of quadratic forms under a mild moment condition. 
\begin{proposition}[Adapted from Lemma A.3 in \citep{el2010spectrum}]
	\label{prop:quad-concentration}
	Let $\{Z_i\}_{i=1}^n$ be i.i.d. random vectors in $\mathbb{R}^d$, whose entries are i.i.d., mean $0$, variance $1$ and $|Z_i(k)| \leq C \cdot d^{\frac{2}{8+m}}$.  
	For any $\Sigma$ under the assumption (A.1), with $\theta = \frac{1}{2} - \frac{2}{8+m}$, we have with probability at least $1 - d^{-2}$,
	\begin{align}
		\max_{i,j} \left| \frac{Z_i^* \Sigma Z_j}{d} - \delta_{ij} \frac{\tr(\Sigma)}{d} \right| \leq d^{-\theta}  \log^{0.51} d,
	\end{align}
	for $d$ large enough.
\end{proposition}
\begin{proof}
	The proof follows almost exactly as in Lemma A.3 \citep{el2010spectrum}. The only point of clarification is that one can assert
	\begin{align}
		P \left( \max_{i,j} \left| \frac{Z_i^* \Sigma Z_j}{d} - \delta_{ij} \frac{\tr(\Sigma)}{d} \right| > 3 d^{-\theta} (\log d)^{0.51}  \right) \leq 2 n^2 \exp(- c (\log d)^{1.02}),
	\end{align}
	and thus for $d$ large enough, say $c(\log d)^{1.02} \geq 4 \log d + \log 2/c^2$, we have that $2 n^2 \exp(- c (\log d)^{1.02}) \leq d^{-2}$. 
\end{proof}

The following proposition is a non-asymptotic adaptation of Theorem 2.1 in \citep{el2010spectrum}. Our contribution here is only to carefully spell out the terms and emphasize that the error rate can be very slow (this is why \citep{el2010spectrum} only provides a convergence in probability result). 
\begin{proposition}
	\label{prop:Karoui-result}
	Under the assumptions (A.1), (A.2), and (A.4), for $\theta= \frac{1}{2} - \frac{2}{8+m}$, with probability at least $1 - \delta - d^{-2}$, 
	\begin{align}
		\left\| K(X, X) - K^{\rm lin}(X, X) \right\| \leq d^{-\theta} (\delta^{-1/2} + \log^{0.51} d),
	\end{align} 
	for $d$ large enough and $\delta$ small enough.
\end{proposition}
\begin{proof}
	In \citep{el2010spectrum}, the approximation error can be decomposed into first-order term $E_1$ (diagonal approximation), second-order off-diagonal term $E_2$, and third-order off-diagonal approximation $E_3$,
	\begin{align}
		K(X, X) - K^{\rm lin}(X, X) := E_1 + E_2 + E_3 
	\end{align}
	where $\| E_1 \| \leq \frac{1}{2} d^{-\theta}  \log^{0.51} d$, $\| E_3 \| \leq \frac{1}{2} d^{-\theta}  \log^{0.51} d$ with probability at least $1 - d^{-2}$. However, for $E_2$ only convergence in probability is obtained. On Page 19 in \citep{el2010spectrum}, the last line reads
	\begin{align}
		\E \tr(E_2^4) \leq C d^{-4\theta},
	\end{align}
	therefore by Chebyshev bound, we have
	\begin{align}
		P\left( \| E_2 \| \geq  d^{-\theta} \delta^{-1/2} \right) \leq P( \tr(E_2^4) \geq d^{-4\theta} \delta^{-2}) \leq \frac{\E \tr(E_2^4)}{d^{-4\theta} \delta^{-2}}  \leq C \delta^2 \leq \delta.
	\end{align}
	
\end{proof}

\section{Lemmas and Corollaries}

\begin{lemma}[Gaussian case]
	\label{lem:1byn-gaussian}
	Under the assumptions (A.1), (A.4), and that $x_i \sim \mathcal{N}(0, \Sigma_d)$ i.i.d. Then with probability at least $1 - d^{-2}$ with respect to a draw of $X$, 
	\begin{align}
		\E_\mu \| K(\bx, X) - K^{\rm lin} (\bx, X) \|^2
		& \leq d^{-1} \log^{4.1} d,
	\end{align}
	with $d$ large enough.
\end{lemma}

\begin{proof}
	Start with entry-wise Taylor expansion for the smooth kernel, 
	\begin{align*}
		K(\bx, x_j) - K^{\rm lin} (\bx, x_j) &= \frac{h''(\xi_j)}{2}  \left(\frac{\bx^* x_j}{d} \right)^2 \leq \frac{M}{2}  \left(\frac{\bx^* x_j}{d} \right)^2.
	\end{align*}
	Conditionally on $x_j, 1\leq j \leq n$, with probability at least $1-2\exp(-t^2/2)$ on $\bx$ drawn from $\mathcal{N}(0, \Sigma)$,
	\begin{align*}
		\left| \frac{\bx^* x_j}{d} \right| = \left| \frac{\langle \Sigma^{1/2} x_j, \Sigma^{-1/2} \bx \rangle}{d} \right| \leq \| \Sigma \| \frac{\| \Sigma^{-1/2} x_j \|}{\sqrt{d}} \frac{t + \log^{0.51} d}{\sqrt{d}}, ~\forall j.
	\end{align*}
	Using standard $\chi^2$ concentration bound, we know that with probability at least $1 - d^{-2}$ on $X$
	\begin{align}
		\max_j \frac{\| \Sigma^{-1/2} x_j \|^2}{d} \leq 1 + \frac{\log^{0.51} d}{\sqrt{d}}.
	\end{align}
	Therefore with probability at least $1-2\exp(-t^2/2)$ on $\bx \sim \mu$, conditionally on $x_j, 1\leq j\leq n$, we have
	\begin{align*}
		\| K(\bx, X) - K^{\rm lin} (\bx, X) \| &\leq C_1 \sqrt{d}  \max_j \left( \frac{\bx^* x_j}{d} \right)^2   \\
		& \leq C_2 \sqrt{d} \cdot \max_j \frac{\|  \Sigma^{-1/2} x_j \|^2}{d}  \cdot  d^{-1} (t^2 + \log^{1.02} d) \\
		& \leq C_3 \max_j \frac{\|  \Sigma^{-1/2} x_j \|^2}{d}  \cdot  d^{-1/2} (t^2 + \log^{1.02} d) \\
		& \leq C_4 \cdot  d^{-1/2} (t^2 + \log^{1.02} d).
	\end{align*}
	Define $Z = \| K(\bx, X) - K^{\rm lin} (\bx, X) \|$. The above says that, conditioned on $X$, $P\left(Z \geq C \cdot  d^{-1/2} (s + \log^{1.02} d)\right) \leq 2\exp(-s/2)$ for all $s>0$. 
	Therefore, by defining change of variables $z = C \cdot  d^{-1/2} (s + \log^{1.02} d)$,
	\begin{align*}
		\E_{\bx \sim \mu} \| K(\bx, X) - K^{\rm lin} (\bx, X) \|^2 &=  \E_{\bx \sim \mu} [Z^2] = \int_{\mathbb{R}_+} 2 z \cdot P(Z \geq z) dz \\
	 	&\leq C \int_{\mathbb{R}_+} d^{-1} (s+ \log^{1.02}d) \cdot \exp(-s/2) ds  \\
		& \leq d^{-1} \log^{4.1} d
	\end{align*}
	with probability at least $1-d^{-2}$ on $X$, for $d$ large enough. 
\end{proof}

\begin{lemma}[Weak moment case]
	\label{lemma:1byn}
	Under the assumptions (A.1), (A.2), and (A.4), for $\theta = \frac{1}{2} - \frac{2}{8+m}$, we have with probability at least $1 - d^{-2}$ with respect to the draw of $X$, for $d$ large enough,
	\begin{align}
		\E_\mu \| K(\bx, X) - K^{\rm lin} (\bx, X) \|^2
		& \leq d^{-(4 \theta - 1)} \log^{4.1} d.
	\end{align}
\end{lemma}
\begin{proof}
	We start with entry-wise Taylor expansion for the smooth kernel, 
	\begin{align*}
		K(x, x_j) - K^{\rm lin} (x, x_j) &= \frac{h''(\xi_j)}{2}  \left(\frac{x^* x_j}{d} \right)^2 \leq \frac{M}{2}  \left(\frac{\bx^* x_j}{d} \right)^2.
	\end{align*}
	Conditionally on $X_j, 1\leq j\leq n$, by Bernstein's inequality \cite[p. 38]{boucheron2013concentration}, with probability at least $1-\exp(-t)$ on $\bx$, for all $j \in [n]$
	\begin{align*}
		\left| \frac{\bx^* x_j}{d} \right| & = \left| \frac{\langle \Sigma^{1/2} x_j, \Sigma^{-1/2} \bx \rangle}{d} \right| \\
		&\leq \sqrt{\frac{2\| \Sigma^{1/2} x_j \|^2}{d}} \frac{\sqrt{t} + \log^{0.51} d}{\sqrt{d}} + \frac{1}{3} \frac{\| \Sigma^{1/2} x_j \|_\infty d^{\frac{2}{8+m}} (t+\log^{1.02} d)}{d},\\
		& \leq   \sqrt{\frac{2\|\Sigma^{1/2} x_j \|^2}{d}} \frac{\sqrt{t} + \log^{0.51} d}{\sqrt{d}} + \frac{1}{3} \frac{\| \Sigma^{1/2} x_j \| d^{\frac{2}{8+m}} (t+\log^{1.02} d)}{d}\\
		& = \frac{\sqrt{2} \| \Sigma^{1/2} x_j \|}{\sqrt{d}}\frac{\sqrt{t} + \log^{0.51} d}{\sqrt{d}} + \frac{1}{3} \frac{\| \Sigma^{1/2} x_j \|}{\sqrt{d}} d^{\frac{2}{8+m} - \frac{1}{2}} (t+\log^{1.02} d) \\
		&= \frac{\| \Sigma^{1/2} x_j \|}{\sqrt{d}} \left( \sqrt{2} d^{-1/2} (\sqrt{t} + \log^{0.51} d) + \frac{1}{3} d^{-\theta} (t+\log^{1.02} d) \right).
	\end{align*}
	Here the second line uses the fact that $\max_k \left|[\Sigma^{1/2} x_j](k) \cdot [\Sigma^{-1/2} \bx](k) \right| \leq \| \Sigma^{1/2} x_j \|_\infty d^{\frac{2}{8+m}}$ due to the assumption~(A.2) for each entry $[\Sigma^{-1/2} \bx](k) \leq C d^{\frac{2}{8+m}}$.
	Applying Proposition~\ref{prop:quad-concentration} with the matrix taken to be identity, for all $j$, with probability at least $1-d^{-2}$ on $X$
	\begin{align*}
		\max_j \frac{\|\Sigma^{1/2} x_j \|^2}{d} \leq \| \Sigma\| \max_j \frac{\|\Sigma^{-1/2} x_j \|^2}{d} \leq C( 1 + d^{-\theta} \log^{0.51} d ).
	\end{align*}
	Therefore with probability at least $1-\exp(-t)$ with respect to $\bx \sim \mu$, conditionally on $x_j, 1\leq j\leq n$ 
	\begin{align*}
		\| K(\bx, X) - K^{\rm lin} (\bx, X) \| &\leq C_1 \sqrt{d}  \max_j \left( \frac{\bx^* x_j}{d} \right)^2  \\
		& \leq C_2 \sqrt{d} \max_j \frac{\| \Sigma^{1/2} x_j \|^2}{d}  \left(  d^{-1} (t + \log^{1.02} d) +   d^{-2\theta} (t + \log^{1.02} d)^2 \right)\\
		\text{(recall $\theta \leq 1/2$)} \quad & \leq C_3 \sqrt{d} \max_j \frac{\| \Sigma^{1/2} x_j \|^2}{d} \left( d^{-2\theta} (t^2 + \log^{2.04} d)  \right) \\
		& \leq C_4  \cdot  d^{-2\theta+1/2} (t^2 + \log^{2.04} d).
	\end{align*}
	Define $Z = \| K(\bx, X) - K^{\rm lin} (\bx, X) \|$. The above says that, conditioned on $X$, $P\left(Z \geq C \cdot  d^{-2\theta+ 1/2} (s^2 + \log^{2.04} d)\right) \leq 2\exp(-s)$ for all $s>0$. 
	Therefore, by change of variables $z = C \cdot  d^{-2\theta+ 1/2} (s^2 + \log^{2.04} d)$, the expectation satisfies
	\begin{align*}
		\E_{\bx \sim\mu} \| K(\bx, X) - K^{\rm lin} (\bx, X) \|^2 &=  \E_{\bx \sim \mu} [Z^2] = \int_{\mathbb{R}_+} 2 z \cdot P(Z \geq z) dz \\
		&\leq  C \int_{\mathbb{R}_+} d^{-4\theta+1} (s^2+ \log^{2.04}d) \cdot \exp(-s) 2s ds   \\
		&\leq C \int_{\mathbb{R}_+} d^{-4\theta+1} s^3 \log^{2.04}d \cdot \exp(-s) ds \\
		& \leq d^{-4\theta+1} \log^{4.1} d
	\end{align*}
	with probability at least $1-d^{-2}$ on $X$, for $d$ large enough.

\end{proof}

\begin{lemma}
	For $\| \Sigma \|_{\rm op} \leq 1$, we have
	\begin{align}
		\label{pf:example-nonparam}
		&\left\| \left[ d {\bf r} \Sigma^{-1} + Z^* Z \right]^{-1} Z^*Z \left[ d{\bf r} \Sigma^{-1} + Z^* Z \right]^{-1} \right\|_{\rm op} \\
		&\leq \left\| \left[ d {\bf r} I +  Z^* Z \right]^{-1} Z^*Z \left[ d {\bf r} I + Z^* Z \right]^{-1} \right\|_{\rm op}. \nonumber
	\end{align}
\end{lemma}
\begin{proof}
	If $A \succeq B$, then $\| A^{-1} C \|_{\rm op} \leq \| B^{-1} C \|_{\rm op} $. Since $\Sigma^{-1} \succeq I$, it holds that
	\begin{align*}
		\left\| \left[ d {\bf r} \Sigma^{-1} +  Z^* Z \right]^{-1} \left[ Z^*Z \right]^{1/2} \right\|_{\rm op} \leq \left\| \left[ d {\bf r} I +  Z^* Z \right]^{-1} \left[ Z^*Z \right]^{1/2} \right\|_{\rm op}.
	\end{align*}
	
\end{proof}

\begin{proof}[Proof of Corollary~\ref{coro:general-n-d}]
	For the variance part, with $Z = X \Sigma_d^{-1/2}$, we have
	\begin{align*}
			{\bf V} &\precsim  \sum_{j} \frac{\lambda_j\left( X^* X  \right)}{ \left[ d {\bf r}   + \lambda_j\left( X^* X   \right) \right]^2}\\
			&=  \tr\left( \Sigma_d^{-1} \left[ d {\bf r} \Sigma_d^{-1} + Z^* Z \right]^{-1} Z^*Z \left[ d {\bf r} \Sigma_d^{-1} + Z^* Z \right]^{-1}  \right) \\
			&\leq  \tr( \Sigma_d^{-1} ) \left\| \left[ d {\bf r} \Sigma_d^{-1} + Z^* Z \right]^{-1} Z^*Z \left[ d {\bf r} \Sigma_d^{-1} + Z^* Z \right]^{-1} \right\|_{\rm op} \\
			&\leq  \tr( \Sigma_d^{-1} ) \left\| \left[ d {\bf r} I +  Z^* Z \right]^{-1} Z^*Z \left[ d {\bf r} I + Z^* Z \right]^{-1} \right\|_{\rm op},
	\end{align*}
	where the last step uses \eqref{pf:example-nonparam}.
	Therefore, using standard random matrix theory, one can further upper bound the above equation by
	\begin{align}
		\label{eq:n>d}
		& \left\| \left[ d {\bf r} I +  Z^* Z \right]^{-1} Z^*Z \left[ d {\bf r} I + Z^* Z \right]^{-1} \right\|_{\rm op} \\
		& = \max_j \frac{\lambda_j\left( Z^* Z \right)}{ \left( d {\bf r} + \lambda_j\left( Z^* Z \right) \right)^2} \leq \frac{1}{n} \frac{(1-\sqrt{d/n})^2}{\left( \frac{d}{n} {\bf r} + (1-\sqrt{d/n})^2 \right)^2} \leq \frac{2}{n} 
	\end{align}
	for $n\gg d$.
	For the bias part,
	\begin{align*}
			{\bf B} &\precsim {\bf r} + \frac{1}{n} \sum_{j=1}^d \lambda_j\left(\frac{X^* X}{d}\right) = {\bf r} + \frac{1}{d} \frac{1}{n} \tr(X^* X) \\
			& =  {\bf r} + \frac{1}{d} \frac{1}{n} \sum_{j=1}^n \| \Sigma^{1/2} z_j \|^2 \nonumber 
			 \precsim {\bf r} + \frac{1}{d} \left( \tr(\Sigma) + \sqrt{\frac{\tr(\Sigma^2)}{n}}  \right) \nonumber
	\end{align*}
	where the last line uses the fact
	\begin{align*}
		\E_{\bf z}\| \Sigma^{1/2} {\bf z} \|^2 = \tr(\Sigma)
	\end{align*}
	and standard $\chi^2$ concentration.
\end{proof}

\begin{proof}[Proof of Corollary~\ref{coro:general-d-n}]
	For the variance bound, we have
	\begin{align*}
		{\bf V} &\leq \frac{n}{d} \frac{1}{4 {\bf r}}, \quad  \text{as $\frac{t}{({\bf r} + t)^2} \leq \frac{1}{4{\bf r}}$ for all $t$}.
	\end{align*}
	For the bias bound, we have
	\begin{align}
		\label{pf:fast-decay}
		{\bf B} &\precsim {\bf r} + \frac{1}{n} \sum_{j=1}^n \lambda_j\left(\frac{X X^* }{d}\right) =  {\bf r} + \frac{1}{d} \frac{1}{n} \sum_{j=1}^n \| \Sigma^{1/2} z_j \|^2 \nonumber \\
		& \precsim {\bf r} + \frac{1}{d} \left( \tr(\Sigma) + \sqrt{\frac{\tr(\Sigma^2)}{n}}  \right) \nonumber
	\end{align}
	by the same argument as in the proof of Corollary~\ref{coro:general-n-d}.
\end{proof}

\begin{lemma}
	\label{lem:symmetrization}
	Let $g(x) \in \mathbb{R}$ that satisfies $\forall g \in \cG$, $|g(x)| \leq M$ for all $x$. Then with probability at least $1-2\delta$, we have for i.i.d. $x_i \sim \mu$
	\begin{align}
		\sup_{g \in \cG} \left| \E g(\bx) - \widehat{\E}_n g(\bx) \right| &\leq \E \sup_{g \in \cG} \left| \E g(\bx) - \widehat{\E}_n g(\bx) \right| + M
		 \sqrt{\frac{\log 1/\delta}{2n}} \\
		 & \leq 2\E \sup_{g \in \cG} \frac{1}{n}\sum_i \epsilon_i g(x_i) + M
		 \sqrt{\frac{\log 1/\delta}{2n}} \\
		 & \leq 2\E_{\epsilon} \sup_{g \in \cG} \frac{1}{n}\sum_i \epsilon_i g(x_i) + 3 M
		 \sqrt{\frac{\log 1/\delta}{2n}}
	\end{align}
	where $\E_{\epsilon}$ denotes the conditional expectation with respect to i.i.d. Rademacher random variables $\epsilon_1,\ldots,\epsilon_n$.
\end{lemma}

\begin{proof}
	The proof is a standard exercise using McDiarmid's inequality and symmetrization. We include here for completeness. See \cite[Theorem 2.21, 2.23 and their corollaries]{mendelson2003few}.
\end{proof}

\newpage
\section{MNIST Result}
\label{sec:mnist}
Here the error is in percentage, so 2.921 corresponds to an error 2.921\%.
{\small
\begin{verbatim}
	Digits pair: [i, j] Error: [Lambda=0   Lambda=0.1 Lambda=1]
	digits pair: [0, 1] error: [0.541      1.006      1.710]
	digits pair: [0, 2] error: [2.921      4.689      7.584]
	digits pair: [0, 3] error: [1.601      3.386      5.841]
	digits pair: [0, 4] error: [1.285      2.610      4.019]
	digits pair: [0, 5] error: [2.567      4.957      8.226]
	digits pair: [0, 6] error: [2.969      5.239      8.359]
	digits pair: [0, 7] error: [1.218      2.808      4.810]
	digits pair: [0, 8] error: [2.541      3.725      5.526]
	digits pair: [0, 9] error: [2.031      3.726      5.482]
	digits pair: [1, 2] error: [2.487      3.699      7.220]
	digits pair: [1, 3] error: [1.644      2.688      4.913]
	digits pair: [1, 4] error: [1.221      2.089      3.552]
	digits pair: [1, 5] error: [1.455      2.860      4.904]
	digits pair: [1, 6] error: [1.615      2.438      3.913]
	digits pair: [1, 7] error: [2.157      3.693      5.689]
	digits pair: [1, 8] error: [2.468      3.571      7.486]
	digits pair: [1, 9] error: [1.441      2.513      3.941]
	digits pair: [2, 3] error: [4.713      7.853     13.253]
	digits pair: [2, 4] error: [2.998      5.602      9.525]
	digits pair: [2, 5] error: [2.711      5.471     10.491]
	digits pair: [2, 6] error: [3.287      5.917     10.519]
	digits pair: [2, 7] error: [4.836      6.930     10.530]
	digits pair: [2, 8] error: [5.080      8.460     13.531]
	digits pair: [2, 9] error: [2.958      5.335      8.763]
	digits pair: [3, 4] error: [1.783      3.847      6.880]
	digits pair: [3, 5] error: [6.822     10.129     17.565]
	digits pair: [3, 6] error: [2.017      3.887      7.088]
	digits pair: [3, 7] error: [3.184      5.486      8.963]
	digits pair: [3, 8] error: [5.345      9.766     16.442]
	digits pair: [3, 9] error: [3.909      6.494     10.330]
	digits pair: [4, 5] error: [2.254      4.871      8.757]
	digits pair: [4, 6] error: [2.878      4.793      7.396]
	digits pair: [4, 7] error: [3.711      7.036     11.015]
	digits pair: [4, 8] error: [3.488      5.615      8.888]
	digits pair: [4, 9] error: [10.199    11.587     18.058]
	digits pair: [5, 6] error: [5.014      7.716     12.682]
	digits pair: [5, 7] error: [2.537      4.683      8.268]
	digits pair: [5, 8] error: [5.868      9.587     16.261]
	digits pair: [5, 9] error: [4.562      6.578     10.935]
	digits pair: [6, 7] error: [1.114      2.864      4.894]
	digits pair: [6, 8] error: [4.102      5.954      9.265]
	digits pair: [6, 9] error: [1.267      2.944      4.935]
	digits pair: [7, 8] error: [3.197      5.623      9.093]
	digits pair: [7, 9] error: [6.598     10.841     17.252]
	digits pair: [8, 9] error: [4.640      7.673     12.070]
\end{verbatim}
}   

\end{document}